\documentclass[10pt]{article}
\usepackage[latin1]{inputenc}
\usepackage{graphicx}
\usepackage{amssymb}
\usepackage{amsfonts}
\usepackage[thmmarks]{ntheorem}
\usepackage[all]{xy}

\def\RR{{\mathbb R}}
\def\CC{{\mathbb C}}
\def\NN{{\mathbb N}}

\def\bra{\langle}
\def\ket{\rangle}
\def\id{{\rm Id}}
\def\tr{\mbox{\rm Tr}}

\def\diag{{\rm diag}}

\def\bea{\begin{eqnarray}}
\def\eea{\end{eqnarray}}

\def\be{\begin{equation}}
\def\ee{\end{equation}}

\newtheorem{theorem}{Theorem}
\newtheorem{definition}{Definition}
\newtheorem{lemma}{Lemma}
\newtheorem{cor}{Corollary}
\newtheorem{propo}{Proposition}
\theoremstyle{nonumberplain}
\theorembodyfont{\normalfont}
\theoremseparator{:}
\theoremsymbol{$\P$}
\newtheorem{demo}{Proof}
\makeatletter
\newcommand{\manuallabel}[2]{\def\@currentlabel{#2}\label{#1}}
\makeatother

\manuallabel{carac1}{3}
\manuallabel{sec26}{2.6}
\manuallabel{sec33}{3.3}
\manuallabel{WRglobal}{3.5}
\manuallabel{eq16}{16}
\manuallabel{propo4}{4}
\manuallabel{lem8}{8}
\manuallabel{lem9}{9}
\manuallabel{causcond}{3.4.3}

\newenvironment{remark}{\small {\bf Remark:}}{}

\begin{document}
\title{On the definition of spacetimes in Noncommutative Geometry: part II}
\author{Fabien Besnard}
\maketitle
\begin{abstract}
In this second part of the paper, we  define  \emph{spectral spacetimes}. They are meant to be the noncommutative generalization of Lorentzian (or antilorentzian, the two cases are worked out) orientable spacetimes of even dimension with a spin structure. There are two main differences with spectral triples: the existence of noncommutative time-orientation  1-forms, and the non-existence of a distinguished    $C^*$-structure on the algebra of ``noncommutative functions''. If a so-called \emph{reconstructibility condition} is met, different, yet isomorphic, $C^*$-structures exist, and an isomorphism is induced by a ``change of observer''.   We define a notion of \emph{stable causality} for  spectral spacetimes which makes a link with previous work by Franco and collaborators. We give an   example of commutative spectral spacetime  which is a Wick rotated version of the spectral triple that must be used in order  to recover the usual notion of distance on a finite graph through Connes' distance formula. We show that this spectral spacetime is stably causal iff the time-orientation 1-form induces no cycle. We also provide a noncommutative example, \emph{the split Dirac structure}, which we study in some details. This structure is defined thanks to a discrete spinor bundle on a finite graph,  a discrete connection on it, and operators at the vertices playing the role of gamma matrices. We give necessary and sufficient conditions for the split Dirac structure to be a spectral spacetime. These includes the natural  analogues of the properties defining a spin connection in the continuous case. However the split Dirac structure is not always reconstructible: we prove that this is the case exactly when there exists a parallel, i.e. covariantly constant, timelike vector field on the graph. Nonreconstructible split Dirac structures furnish examples of spectral spacetimes which are not Wick rotations of usual spectral triples, and are thus genuinely Lorentzian. Moreover we show that the Dirac operator of the split Dirac structure is related to a Lorentzian version of a  discretization of the   Dirac operator introduced earlier by Marcolli and van Suijlekom. 
% which turns out to be related to a Lorentzian version of the discretization of the   Dirac operator introduced by Marcolli and van Suijlekom. We study these examples in some detail, using them in particular to argue against the idea that the algebra of ``noncommutative functions'' should be constrained to be stable under the Krein adjunction or, even, be a $C^*$-algebra. Finally we use the split Dirac structure as a toy-model to show the potential physical applications of the approach. Namely, we show that in dimension 4 the spectral action is extremal on disjoint unions of   graphs with lightlike edges and no oriented cycle.
\end{abstract}

\tableofcontents
\section{Introduction to part II}
In this second part of the paper we put forward the general definition of  the structures we hope to be the noncommutative counterparts of time and space oriented spin Lorentzian/antilorentzian  manifolds of even dimension. We stay at an algebraic level, setting aside the analytical questions. For this reason the examples we will give, with the exception of the classical case of a manifold, will be finite-dimensional.  We deal almost exclusively with the antilorentzian case. We refer to \cite{part1} for a general introduction to this work and notations.

This part of the paper is organized as follows: in section 2 we give the definition \emph{spectral spacetimes}, both in the Lorentzian and antilorentzian signature. We define the \emph{reconstructibility condition}, under which the algebra will have at least one $C^*$-structure. The notion of \emph{stable causality} is also straightforwardly generalized from manifolds to spectral spacetime. We compare spectral spacetimes with previous approaches to Lorentzian noncommutative geometry. Finally we explain how (and when) back and forth Wick rotations between spectral spacetimes and spectral triples can be performed.

In section 3 we focus on finite-dimensional spectral spacetimes.  First we consider the issue of stable causality, and find it inconsistent with the   first-order condition which is  usually postulated in noncommutative geometry. We then give three finite-dimensional examples of spectral spacetimes, each being a generalization of the previous one. The first example is the simplest possible, with algebra $\CC^2$. We find essentially only one possibility, which has KO dimension 2. We note that it contrasts with spectral triples since two families of spectral triples exists with algebra $\CC^2$, one of KO dimension 0 and one of KO dimension 6. We point out that the spectral spacetime we found is the Wick rotation of the 0-dimensional spectral triple while the 6-dimensional one cannot be Wick rotated. The second example is the canonical antilorentzian spectral spacetime over a positively weighted directed graph. We prove that it is stably causal iff the orientation on the graph is acyclic. We note that this example is not covered by the various   approaches reviewed in subsection \ref{sec44}. The last example is a noncommutative generalization of the previous one: the split Dirac structure. We characterize the cases where this structure gives a reconstructible spectral spacetime. We use this to find examples of spectral spacetimes which are not reconstructible. In particular they are not Wick rotations of spectral triples. We also perform a preliminary study of the stable causality condition, and show the relation between the split Dirac structure and the discretization of the usual Dirac operator.
 %However, we note that this example is far from being far-fetched, as it is precisely related with the discretization of the usual Dirac operator.

Section \ref{conclusion} contains a conclusion, outlooks and acknowledgements.

\section{Lorentzian and anti-Lorentzian spectral spacetimes}\label{sec5}
\subsection{Preliminary: noncommutative 1-forms}
In order to define a noncommutative generalization of (space and time-oriented, spin, even dimensional) antilorentzian spacetimes we  will take inspiration from \cite{part1}, theorem \ref{carac1}. Hence we need to use an object $\beta$ which plays the role of a differential 1-form. The definition of noncommutative forms over a spectral triple can be found in \cite{redbook}, chapter VI. However we do not work with self-adjoint operators on Hilbert space, so strictly speaking this theory does not apply, although the modifications we need to make are very mild. We just need  noncommutative 1-forms in the present paper, hence we only quote the relevant definitions. Some more details on general noncommutative forms in a semi-riemannian context can be found in \cite{bbb}.

Let ${\cal A}$ be a unital algebra. The \emph{universal differential graded algebra} over ${\cal A}$, written $\Omega{\cal A}$, is generated by the elements of ${\cal A}$, which are given degree $0$, and symbols $da$, for $a\in {\cal A}$, which are given degree 1, submitted to the relations :

\bea
d(ab)&=&(da)b+adb\cr
d(1)&=&0\label{leibrule}
\eea

In particular the elements of degree one can be reduced to the form $\sum_i a_i db_i$. Their space is denoted by $\Omega^1{\cal A}$ and is  an ${\cal A}-{\cal A}$ bimodule since it is stable by multiplication on both sides by elements of degree $0$ thanks to (\ref{leibrule}). The map $d : a\mapsto da$ can be uniquely extended to $\Omega{\cal A}$, turning it into a differential algebra, by requiring  the graded Leibniz rule and $d^2=0$ to hold.

Now let $\pi$ be a representation of ${\cal A}$ in a space $K$, and $D$ an operator on $K$. This representation extends as a representation $\tilde \pi_D$ of $\Omega{\cal A}$ by the rules $\tilde\pi_D(a)=\pi(a)$, $\tilde\pi_D(da)=[D,\pi(a)]$, for all $a\in {\cal A}$, thanks to the fact that $[D,.]$ is a derivation. The \emph{space of noncommutative 1-forms over $({\cal A}, \pi, D)$} is defined to be $\Omega^1_D({\cal A},\pi):=\tilde\pi_D(\Omega^1{\cal A})$. It is a  sub-${\cal A}-{\cal A}$ bimodule of $End(K)$ since $\Omega^1({\cal A})$ is a sub-${\cal A}-{\cal A}$ bimodule of $\Omega {\cal A}$ and $\tilde\pi_D : \Omega{\cal A}\rightarrow End(K)$ is an algebra morphism.  The noncommutative 1-forms over $({\cal A}, \pi, D)$ are thus elements of the form 
$$\sum_i\pi(a_i)[D,\pi(b_i)]$$

for some $a_i,b_i\in{\cal A}$.  

\subsection{Definition of spectral spacetimes}

%Recall that a fundamental symmetry $\beta$ satisfies $\beta^2=1$, $\beta^\times=\beta$, $\bra .,.\ket_\beta>0$. It is obviously self-adjoint relatively to $\bra .,.\ket_\beta$.

%A global time function $a$ will define at each point of spacetime a preferred observer and a $C^*$-algebra $A_\beta$ of observables associated to this observer, where $\beta=i\gamma(da)$. We ask that they are all isomorphic to a model $C^*$-algebra $A$. Hence they will define isomorphic spaces of states.

We deal first with the antilorentzian case. We will say a word about the required modifications in the Lorentzian case at the end of this section.

Since we set aside analytical properties things are pretty straightforward. Nevertheless, an interesting difficulty appears when we replace the commutative algebra ${\cal C}(M)$ by a noncommutative one: the $C^*$-involution becomes dependent on the positive time-orientation form. Indeed, the adjoint of an operator $A$ on the Krein  space $(K,(.,.))$ with respect with the scalar product $\bra .,.\ket_\beta=(.,\beta^{-1}.)$ is $A^{*_\beta}=\beta A^\times \beta^{-1}$ where $\times$ is the Krein adjoint. In the classical case, the  Krein adjoint of the multiplication operator by a function $f$ on $M$ is the multiplication operator by $\bar f$, and as such commutes with the 1-form $\beta$, thus  $\pi(f)^\times=\pi(f)^{*_\beta}=\pi(\bar f)$. There would be no reason to suppose this to be true in the noncommutative case. On the contrary, it should be recovered as a particular case when the first-order condition holds and the algebra is commutative. The necessity for the $C^*$-involution to become dependent on $\beta$ can also be understood in physical terms: the 1-form $\beta$ is dually equivalent to a non-vanishing timelike vector fields. Since a vector field is integrable, this is also equivalent to a foliation of spacetime by timelike curves, which can be seen as defining a congruence of observers. It is physically sound that the scalar product on the space of states depends on this congruence of observers. This is what happens for instance in relativistic quantum mechanics. One can also consult \cite{rov1} to see how this observation applies in the context of loop quantum gravity. 

Hence we see that we need to abandon the $C^*$-paradigm: we cannot take a particular $C^*$-algebra as a basic object in the definition we are seeking. We will  discuss what will replace it after we give the definition.

% It would be tempting to take a Krein-$C^*$-algebra instead. We refrain from doing so for two reasons: the first is that it would require the image of  algebra in $End(K)$ to stable under the Krein adjoint, and we will see in section 6.4 that there is a very natural example in which this is not the case, the second is that it would put an apriori restriction on the possible time-orientation forms, which is not natural.

\begin{definition}\label{antilorsst} An even  antilorentzian spectral spacetime of $KO$-dimension $k$ mod $8$ is given by a multiplet $S=({\cal A}, K, (.,.),\pi,D,J,\chi)$  where:
\begin{enumerate}
\item the couple  $(K,(.,.))$ is a Krein space, with adjunction denoted by $\times$,
\item ${\cal A}$ is an algebra and $\pi$ is  a   representation of it on $K$,
\item the ``Dirac operator'' $D$ on $K$ is such that $D^\times=D$,
\item the  ``chirality operator'' $\chi\in B(K)$ is such that $\chi^2=1$, $[\pi(a),\chi]=0$ for all $a\in A$, $\chi D=-D\chi$ and $\chi^\times=-\chi$,
%\item A universal one-form $\delta\in\Omega^1({\cal A})$ such that $\pi_D(\delta)=\beta$.
\item the ``charge conjugation'' $J$  is an antilinear operator on $K$   which is required to satisfy:
$$J^2=\epsilon,\quad JD=DJ,\quad J\chi=\epsilon''\chi J,\quad J^\times J=\kappa$$
with $\epsilon,\epsilon'',\kappa=\pm 1$, according to  table \ref{KOantilor}.
\end{enumerate}
and such that there exists an element $\beta\in \Omega^1_D({\cal A},\pi)$, which satisfies the following properties:
\begin{enumerate}
\item it is Krein self-adjoint,
\item it is  imaginary ( $J\beta J^{-1}=-\beta$),
\item the hermitian form    $\bra .,.\ket_\beta:=( .,\beta^{-1}.)$ is  positive definite,
\end{enumerate}
\end{definition}

An element $\beta$ satisfying $1,2,3$ above is called a \emph{positive  time-orientation 1-form}.

Of course the name ``spectral spacetime'' is a bit sloppy, since our definition aims at corresponding to the classical notion of an even dimensional spin antilorentzian manifold with a space and time orientation (``physics-friendly'' spacetimes). But we adopt it for the sake of simplicity. Let us also remark that  since only 4 out of the 8 possibilities of triple of signs appear in the KO-dimension table of signs, the possibility of finding a structure satisfying all the axioms of a spectral spacetime except for these signs is not excluded.   We haven't found any such example yet.

\begin{table}
\begin{center}
\begin{tabular}{ccccc}
metric dim $[8]$ & 0 & 2 & 4 & 6 \\
KO dim $[8]$ & 2 & 0 & 6 & 4 \\
$\epsilon$ & -1 & 1 & 1 & -1 \\
$\epsilon''=\tilde\epsilon''$ & -1 & 1 & -1 & 1 \\
$\tilde\epsilon=\epsilon''\epsilon$ & 1 & 1 & -1 & -1 \\
%$\tilde\epsilon'$ & -1 & -1 & -1 & -1 \\
$\kappa$ & -1 & -1 & -1 & -1 \\
$\tilde\kappa$ & 1 & 1 & 1 & 1 
\end{tabular}
\caption[smallcaption]{Antilorentzian table of  KO signs.}
\label{KOantilor}
\end{center}
\end{table}

In order to replace the $C^*$-structure on ${\cal A}$ we introduce the following notion:

\begin{definition} If $S$ is a spectral spacetime such that $\pi$ is faithful and there exists a positive time-orientation form such that 

\be
 \pi({\cal A})^{*_\beta}=\pi({\cal A})\label{cond4def}
\ee
we say that $S$ is \emph{reconstructible}.
\end{definition}

If $S$ is reconstructible then thanks to the injectivity of $\pi$ we can transport $*_\beta$ on ${\cal A}$ and this  endows it with a $C^*$-algebra structure such that $\pi$ is a $C^*$-morphism. While different time-orientation forms $\beta$ and $\beta'$  induce different  $C^*$-structures on ${\cal A}$, these are canonically isomorphic. Indeed, $Ad_{\beta'\beta^{-1}}$  stabilizes $\pi({\cal A})$ and sends $\pi(a)^{*_\beta}$ to $\pi(a)^{*_{\beta}'}$. Thus, these $C^*$-structures allow the reconstruction of homeomorphic (and hopefully, under convenient assumptions, \emph{diffeomorphic}) structure spaces. This explains the name of this condition which would certainly be needed for a reconstruction theorem. \emph{In the rest of the paper we will always assume reconstructibility except when specified otherwise.}

We now make a few more comments about the definition, and introduce some terminology.

\begin{itemize}
\item It is easy to see that if $\beta$ is a positive time-orientation form then the other positive time-orientation forms can be expressed as $\omega=\rho\beta$ where $\rho$ is a positive operator on the Hilbert space $(K,\bra .,.\ket_\beta)$ which commutes with  $J$. Moreover the reconstructibility condition will hold  for $\omega$ iff $\rho\pi({\cal A})\rho^{-1}\subset\pi({\cal A})$.

%The set of all such 1-forms is by definition the positive time-orientation of the spectral spacetime and is denoted by $Or_+^t(S)$.
\item If we decompose $K$ as a sum $K_L\oplus K_R$ of $-1/+1$-eigenspaces of $\chi$, respectively, we   see from $\chi^\times=-\chi$ that $K_L$ and $K_R$ are self-orthogonal for the Krein product. On the other hand, since every 1-form anticommutes with $\chi$, we see that $K_L$ and $K_R$ are orthogonal to each other for the scalar product $\bra .,.\ket_\beta$ defined by any positive time-orientation form.
\item As in \cite{part1}, section \ref{sec26} it will be convenient to say that $\beta$ is \emph{Krein-positive} whenever $(.,\beta .)$ is positive definite. Since $\beta$ is Krein-positive iff $\beta^{-1}$ is, and Krein-positivity implies Krein self-adjointness, we could have just said in   definition \ref{antilorsst} that a positive orientation-form is Krein-positive and imaginary. Note however that the reconstructibility  condition  is not symmetrical under the exchange $\beta\leftrightarrow \beta^{-1}$. Moreover $\beta^{-1}$ is not guaranteed to be a 1-form.

\item If in addition to reconstructibility we ask that $\pi({\cal A})^\times=\pi({\cal A})$, in which case we say that $\pi({\cal A})$ is \emph{a $\times$-subalgebra of $B(K)$}, then ${\cal A}$ will be equipped with an involution, still denoted $\times$, and will be a Krein $C^*$-algebra. It would be tempting to ask this property as an additional axiom, since it is satisfied in the classical case. However, we will see that it is not satisfied in the discrete commutative case of a positively weighted graph (see section \ref{sec64}). Since the continuous and the discrete must be on the same footing in noncommutative geometry, it would not be appropriate to ask for $\pi({\cal A})$ to be a $\times$-subalgebra in general.
%\item Let us observe, as a side remark, that the characterization of the signature   uses every object in the list $S=({\cal A}, K, (.,.),\pi,D,J,\chi)$.
\end{itemize}

\begin{remark} We could have thought to require the property $\omega^\times=-J\omega J^{-1}$ for general 1-forms in order to mimick the $c$-compatibility of the representation (this is equivalent to real 1-forms being Krein-self-adjoint in the classical case), but this would not make sense in general (under the order 1 condition it would imply commutation between 1-forms and the algebra). Instead the $c$-compatibility is taken care of by the requirement $D=D^\times$ as explained in \cite{part1}, section \ref{sec33}.
\end{remark}
\smallbreak
% We might require additional properties like:

% \be
% \beta^2=1\label{normalization}
% \ee

% and/or

% \be
% J\beta J^{-1}=-\beta\label{imaginarity}
% \ee

% When (\ref{normalization}) is satisfied we say that $\beta$ is normalized, and when (\ref{imaginarity}) is satisfied we say that it is imaginary. These two properties exactly correspond to the one we imposed on the Pin group element $g$ in proposition \ref{generalWR}.

Let us speak now about the \emph{Lorentzian spectral spacetimes}. They are defined exactly as the antilorentzian ones, except for the following points:

\begin{enumerate}
\item The hermitian form $\bra .,.\ket_\beta$ is defined to be $(.,\beta^{-1}\chi.)$.
\item The signs $\epsilon,\epsilon'',\kappa$ are given by table \ref{KOlor}.
\end{enumerate}

It is immediate to check that $\beta^{-1}\chi$ is Krein self-adjoint iff $\beta$ is. 

\begin{table}
\begin{center}
\begin{tabular}{ccccc}
metric dime  $[8]$& 0 & 2 & 4 & 6 \\
KO dim $[8]$ & 6 & 0 & 2 & 4 \\
$\epsilon$ & 1 & 1 & -1 & -1 \\
%$\epsilon'$ & 1 & 1 & 1 & 1 \\
$\epsilon''=\tilde\epsilon''$ & -1 & 1 & -1 & 1 \\
$\tilde\epsilon=\epsilon''\epsilon$ & -1 & 1 & 1 & -1 \\
%$\tilde\epsilon'$ & -1 & -1 & -1 & -1 \\
$\kappa$ & 1 & -1 & 1 & -1 \\
$\tilde\kappa$ & -1 & 1 & -1 & 1 
\end{tabular}
\caption[smallcaption]{Lorentzian  table of KO signs.}
\label{KOlor}
\end{center}
\end{table}

In both the Lorentzian and antilorentzian cases we might require additionally that 

\be
\beta^2=1\label{normalization}
\ee

When (\ref{normalization}) is satisfied we say that $\beta$ is \emph{normalized}. If there exists a normalized positive time-orientation on a spectral spacetime we say that it is \emph{normalizable}.

%Before proceeding let us come back to an issue raised earlier. Since now $\beta$ has no reason to be proportional to $\beta^{-1}$ we could insist on defining $\bra .,.\ket_\beta$ by $(\beta^{-1}.,.)$ instead of $(\beta .,.)$. However this would not change anything: $\beta$ is a self-adjoint operator for $\bra .,.\ket_\beta$, so that $\beta^{-2}$ is a positive operator, hence $\bra \beta^{-2}\psi,\psi\ket_\beta=(\beta^{-1}\psi,\psi)$ is always positive. 

% \begin{lemma} $\pi^o$ is a $*_\beta$-representation of ${\cal A}$ iff $(\beta J)^2$ commutes with $\pi({\cal A})$. 
% \end{lemma}
% \begin{demo}
% We have 
% \bea
% \pi^o(a^{*_\beta})&=&J\pi(a^{*_\beta})^\times J^{-1}\cr
% &=&J(\pi(a)^{*_\beta})^\times J^{-1}\cr
% &=&J(\beta^{-1}\pi(a)^\times\beta)^\times J^{-1}\cr
% &=&J\beta\pi(a)\beta^{-1}J^{-1}
% \eea
% and
% \bea 
% \pi^o(a)^{*_\beta}&=&\beta^{-1}(J\pi(a)^\times J^{-1})^\times \beta\cr
% &=&\beta^{-1}(J^\times)^{-1}\pi(a)J^\times\beta\cr
% &=& \kappa^2\beta^{-1}J\pi(a) J^{-1}\beta
% \eea

% Hence $\pi^o$ is a $*_\beta$-representation iff $\pi(a)$ commutes with $J^{-1}\beta J\beta$ for all $a$. Since $J^{-1}=\epsilon J$ the lemma follows.
% \end{demo}

Finally, by analogy with the classical case (see subsection \ref{causcond} in \cite{part1}) we  define  \emph{stably causal} spectral spacetimes.

\begin{definition} An antilorentzian/Lorentzian spectral space is said to be \emph{stably causal} iff there exists a time-orientation form on it which is exact.
\end{definition}

Recall that the real even spectral triples $S_1=({\cal A},H_1,\pi_1,D_1,J_1,\chi_1)$ and $S_2=({\cal A},H_2,\pi_2,D_2,J_2,\chi_2)$ over the same algebra ${\cal A}$ are called unitarily equivalent iff there exists a unitary operator $U: H_1\rightarrow H_2$ which intertwines $\pi_1$ and $\pi_2$, and such that $UD_1U^*=D_2$, $UJ_1U^*=J_2$ and $U\chi_1U^*=\chi_2$. 

We will define \emph{unitary equivalence} of real even spectral spacetimes in a way which generalizes straightforwardly unitary equivalence of spectral triples. This notion exactly corresponds in the classical case to the action of diffeomorphisms\footnote{preserving the spin structure as well as space and time orientation} of a manifold by pullback.

\begin{definition}  The real even spectral spacetimes $S_1=({\cal A},K_1,\pi_1,D_1,J_1,\chi_1)$ and $S_2=({\cal A},K_2,\pi_2,D_2,J_2,\chi_2)$ over the same algebra ${\cal A}$ are called \emph{unitarily equivalent} iff there exists a Krein-unitary operator $U: K_1\rightarrow K_2$ such that

$$U\pi_1(.)U^{-1}=\pi_2(.) , UD_1U^{-1}=D_2, UJ_1U^{-1}=J_2, U\chi_1U^{-1}=\chi_2$$

and which sends  positive timelike forms to positive timelike forms by its adjoint action.
\end{definition}

%\subsection{Discussion}\label{sec43}

% \subsubsection{$C^*$-algebra structure on ${\cal A}$}

% Most of the time the representation $\pi$ is faithful, in which case ${\cal A}$ can be directly seen as a subalgebra of $B(K)$. The requirement that $\pi$ is a $*_\beta$-morphism then boils down to $\beta^{-1}{\cal A}^\times \beta\subset{\cal A}$. If in addition ${\cal A}^\times={\cal A}$, then ${\cal A}$ naturally becomes a Krein $C^*$-algebra.  

% NO, WRONG : THE REQUIREMENT MADE HERE IS WEAKER, THE C* STRUCTURE MIGHT DEPEND ON BETA WHERAS IN THE DEF IT DOES NOT. (There is one fixed C*structure)

\subsection{Conditions of order $0$ and $1$}\label{sec521}

In this section we explore the generalization to the antilorentzian setting of the conditions of order 0 and 1. We will see that they put severe constraints on the orientation forms. 

Let us define the right representation $\pi^o : {\cal A}\rightarrow B(K)$ by  $\pi^o(a)=J\pi(a)^\times J^{-1}$. Note that if $\pi({\cal A})$ is a $\times$-subalgebra, then $\pi^o$ is an involutive algebra morphism ($\pi^o(a^\times)=\pi^o(a)^\times$).

\begin{remark}
Given a time-orientation form $\beta$ satisfying the reconstructibility condition  we could try to define $\pi^o$ thanks to the Hilbert adjoint instead of the Krein adjoint, but then it would not make sense to use $J$. We could use $\beta J$ and define $\pi^o(a)=\beta J\pi(a)^{*_\beta}(\beta J)^{-1}$. Since $J$ anticommutes with $\beta$ this definition agree with the previous one, and in particular does not depend on $\beta$.
\end{remark}

\begin{propo} Let $\beta$ be a time-orientation form   satisfying the reconstructibility condition. The following properties are equivalent.
\begin{enumerate}
\item $\pi^o$ is a $*_\beta$-representation of ${\cal A}$
\item $\beta^2$ commutes with $\pi({\cal A})$
\item For all $a\in{\cal A}$, $(\pi(a)^{*_\beta})^\times=(\pi(a)^\times)^{*_\beta}$.
%\item The involutions $\times$ and $*_\beta$ commute on ${\cal A}$.
\end{enumerate}
\end{propo}
\begin{demo}
We have 
\bea
\pi^o(a^{*_\beta})&=&J\pi(a^{*_\beta})^\times J^{-1}\cr
&=&J(\pi(a)^{*_\beta})^\times J^{-1}\cr
&=&J(\beta\pi(a)^\times\beta^{-1})^\times J^{-1}\cr
&=&J\beta^{-1}\pi(a)\beta J^{-1}
\eea
and
\bea 
\pi^o(a)^{*_\beta}&=&\beta(J\pi(a)^\times J^{-1})^\times \beta^{-1}\cr
&=&\beta(J^\times)^{-1}\pi(a)J^\times\beta^{-1}\cr
&=& \kappa^2\beta J\pi(a) J^{-1}\beta^{-1}
\eea

Hence $\pi^o$ is a $*_\beta$-representation iff $\pi(a)$ commutes with $\beta J^{-1}\beta J=-\beta^2$ for all $a$. The equivalence between the second and  third properties is immediate. 
\end{demo}

Now the zeroth and first-order conditions read, respectively:

\be
\forall a,b\in{\cal A},\ [\pi(a),\pi^o(b)]=0\label{zeroth}
\ee

and

\be
\forall a,b\in{\cal A},\ [[D,\pi(a)],\pi^o(b)]=0\label{first}
\ee

%The zeroth order condition really makes sense only when $\pi^o$ is $*_\beta$-representation.  

\begin{propo}\label{propo12} If the conditions of order $0$ and $1$ are satisfied, then any time-orientation form $\beta$ commutes with $\pi({\cal A})$.
\end{propo}
\begin{demo}
Let $a\in{\cal A}$. Then
\bea
(J[\beta,\pi(a)]J^{-1})^\times&=&[J\beta J^{-1},J\pi(a)J^{-1}]^\times\cr
&=&-[\beta,J\pi(a)J^{-1}]^\times\cr
&=&[\beta,\pi^o(a)]\cr
%&=&[\beta,\beta^{-1}\pi^o(a)^{*_\beta}\beta]\cr
&=&0
\eea
since by the conditions of order $0$ and $1$ the element $\pi^o(a)$ commutes with any 1-form.
\end{demo}

\subsection{Comparison with previous works}\label{sec44}

The more usual definitions of semi-Riemannian spectral triples found in the literature stem from the seminal work \cite{strohmaier}. In this section we will stress the most important   similarities and differences with our definition of spectral  spacetimes. We quote the relevant definitions with slight changes of notations in order to make the comparison easier. Of course if a spectral spacetime is not reconstructible, a major difference with all the approaches so far is that the algebra is not equipped with a $C^*$-product. Hence in this section we will always take a reconstructible spectral spacetime as a comparison basis.

In \cite{strohmaier} an even \emph{semi-Riemannian spectral triple} is defined to be:

\begin{enumerate}
\item A $*$-algebra ${\cal A}$ and a Krein space $(K,(.,.))$.
\item A faithful representation $\pi$ such that $\pi(a^*)=\pi(a)^\times$.
\item A grading $\chi$ such that $\chi^2=1$, $\chi^\times=\pm\chi$ and $\chi$ commutes with $\pi({\cal A})$.
%\item A fundamental symmetry $\beta$ such that  $\beta^\times=\beta$, $\beta^2=1$ and $\beta$ commutes with $\pi({\cal A})$. %Moreover ${\cal J}$ is required to be either even or odd. When it is off the Krein spectral triple is said to be of Lorentz type.
\item A Dirac operator $D$ which is Krein-self adjoint and anticommutes with $\chi$.
\end{enumerate}

%Moreover a fundamental symmetry $\beta$ is said to be \emph{admissible} if $\beta\chi\beta=\chi^\times$ and $\beta$ commutes with $\pi({\cal A})$. This class of fundamental symmetry is viewed as an analogue of the spatial reflections on a space-oriented manifold.

By rights, spectral  spacetimes should be semi-Riemannian spectral triples, but we see that they are not, since the $*$-structure on the algebra of a reconstructible spectral spacetime corresponds to the Hilbert adjoint, not the Krein adjoint. This is a major difference. In fact we see that a spectral spacetime is a semi-Riemannian spectral triple if and only if there exists a time-orientation form $\beta$ which is normalized and commutes with $\pi({\cal A})$. %In this very special case $\beta$ is an admissible fundamental symmetry.
 
The definition of  \emph{Krein spectral triples of Lorentz type}  given by van den Dungen in  \cite{vdd} corresponds to 
Strohmaier's semi-Riemannian spectral triples with a given  fundamental symmetry $\beta$ which commutes with $\pi({\cal A})$ and is odd.

In the approach of \cite{vddpr}  the Hilbert space is taken to be a basic object and no fundamental symmetry is assumed. The algebra is a $*$-algebra represented on the Hilbert space and the Dirac operator is not self-adjoint. This setting is on the one hand very general, as can be seen already in the finite-dimensional case. In section  4.3.3 of \cite{vddpr} it is observed that a Dirac operator can take the form of a triangular matrix, which is certainly not the case with a spectral spacetime (see section \ref{finitedim} below). On the other hand  fixing a Hilbert structure does not permit to encode in a natural way the covariant family of Hilbert structures which may exist on a reconstructible spectral spacetime.

In all the approaches presented so far there is no axiom which pinpoints the Lorentzian or antilorentzian signature. The sign in  the relation $\chi^\times=\pm \chi$ is sensitive to the signature, but this clearly does not suffice. As far as we are aware the first paper dealing with this matter is \cite{PS}, where \emph{geometric Lorentzian spectral triples} are defined. In short they are semi-Riemannian spectral triples with a fixed fundamental symmetry commuting with the algebra, and a charge conjugation operator. The conditions of order $0$ and $1$ are assumed to hold. The geometric Lorentzian spectral triple is said to be time-oriented when the fundamental symmetry takes on the following form:

\be
\beta=i\sum_i J\pi(a_i^0)J^{-1}\pi(a_i)[D,\pi(b_i)]\label{psform}
\ee

where $a_i,a_i^0,b_i\in{\cal A}$. We note that when $\beta$ is a one-form, the geometric Lorentzian spectral triple is time-oriented. We see then that this definition is very close to the one we have given here\footnote{Our notation $\beta$ comes from this paper. We think this is better that the more usual ${\cal J}$ since it avoids confusion  with the charge conjugation operator.}, with some caveats:

\begin{itemize}
\item There is nothing specifically Lorentzian in the absence of time-orientation.
\item The form (\ref{psform}) for $\beta$ makes sense in a context where the condition of order 1 holds, which is something we want to get away from in view of our study of finite-dimensional examples below.
\item The fundamental symmetry here again commutes with the algebra and is normalized.
\end{itemize}

Unfortunately the full justification of these axioms has been postponed by the authors.

Another attempt at building a specifically Lorentzian framework for noncommutative geometry can be found in the works of Franco and collaborators (\cite{franco}, \cite{francoeckstein}, \cite{francowallet}). We will use the definition of \emph{Lorentzian spectral triples} given in \cite{francoeckstein}. Here $\beta$ is assumed to commute with $\pi({\cal A})$, to be normalized, and moreover  to be of the form $\beta=-N[D,T]$, where $N$ belongs to the algebra and $T$ is a (possibly unbounded) operator on the Hilbert space. They  correspond respectively to lapse and time functions. When $T$ belongs to the algebra we see that $\beta$ is a noncommutative 1-form.

The main differences between our work and other approaches are summarized in table \ref{diff}. Our choice to define $\beta$ to be an imaginary 1-form is motivated by \cite{part1}, theorem \ref{carac1}. Moreover, the postulate that $\beta$ is imaginary has the virtue of making the definition of $\pi^o$ unambiguous (see the remark at the beginning of subsection \ref{sec521}). It also turns out to be a necessary condition in proposition \ref{WAL} below. As a consequence it is needed in order to remove examples which would have bad KO-triple of signs (see example \ref{dim2} below). Of course the most important departure from previous approaches is that we no longer suppose that $\beta$ commutes with the algebra. At the risk of repeating ourselves, we think we are fully justified in doing so by the existence of natural finite-dimensional spectral spaces where this axiom is not satisfied (see  section \ref{sec64}). We admit not to have any justification as strong as these for leaving $\beta$ unnormalized. We can at least invoke the need for generality. We haven't found yet any example where $\beta$ cannot be normalized. Hence it could turn out to be a theorem rather than a postulate that every spectral spacetime is normalizable.

\begin{table}
\begin{minipage}{\textwidth}
 \begin{center}
 \begin{tabular}{ccccc}
 & Stromaier & Paschke-Sitarz & Franco & This work \\
 $\beta$ is fixed & no & yes & yes & no\\
 $\beta$ is normalized & yes & yes & yes & no \\
 $\beta$ commutes with the algebra & yes & yes & yes & no\footnote{Yes when the conditions of order 0 and 1 are satisfied.} \\
 $\beta$ is a 1-form & no & can be\footnote{It is not a 1-form in general, but there are natural cases where it is.} & can be\textsuperscript{\it b} & yes \\
  $\beta$ is imaginary & no & yes & no  & yes \\
 \end{tabular}
 \caption[smallcaption]{Comparison between several  approaches to Lorentzian NCG}\label{diff}
 \end{center}
 \end{minipage}
\end{table}

\subsection{Wick rotations of spectral triples and spectral spacetimes}

Let us first see how we can Wick rotate an antilorentzian spectral spacetime to Euclidean signature. This is closely parallel to \cite{part1}, section \ref{WRglobal}. Note that the reconstructibility condition is assumed in this subsection, and proves to be important.

\begin{propo}\label{ALversEUC} Let $S=({\cal A},K,(.,.),\pi,D,J,\chi)$ be an  antilorentzian spectral spacetime of $KO$-dimension $n$ mod $8$, and let $\beta$ be a positive time orientation 1-form on $S$. Define $J_\beta=\beta J$, $\chi_\beta=-\chi$,  $D_\beta={1-i\over 2}D+{1+i\over 2}\beta D\beta^{-1}$, and $a^{*_\beta}=\pi^{-1}(\pi(a)^{*_\beta})$.

Then $S_\beta=(({\cal A},*_\beta),K,\bra .,.\ket_\beta,\pi,D_\beta,J_\beta,\chi_\beta)$ is a real even  spectral triple with KO-dimension $2-n$   if and only if $\beta$ is normalized. %and imaginary.
%\item If $\beta$ is imaginary then $\pi^o(a)=J_\beta\pi(a)^{*_\beta}J_\beta^{-1}$ for all $a\in {\cal A}$ and $\pi^o$ is $*$-representation of ${\cal A}$ on the Hilbert space $(K,\bra .,.\ket_\beta)$.
\end{propo}
\begin{demo}
   We call $\epsilon,\epsilon'',\kappa$ the KO signs of $S$, and $\epsilon_\beta,\epsilon_\beta'',\kappa_\beta$ the KO signs of the spectral triple we  want to obtain by Wick rotation. We need to have $\epsilon_\beta=-\epsilon$, $\epsilon_\beta''=-\epsilon''$ and $\kappa_\beta=-\kappa=1$. %Let us first prove that it is necessary for $\beta$ to be normalized and imaginary. We call $\epsilon,\epsilon'',\kappa$ the KO signs of the given AL spectral space, and $\epsilon_\beta,\epsilon_\beta'',\kappa_\beta$ the KO signs of the spectral triple we obtained by Wick rotation. We need to have $\epsilon_\beta=-\epsilon$, $\epsilon_\beta''=-\epsilon''$ and $\kappa_\beta=-\kappa=1$. 

%However we have $J_\beta^2=\beta J\beta J=\epsilon \beta J\beta^\times J^{-1}$, since $\beta^\times=\beta$ and $J^{-1}=\epsilon J$. Letting $\beta^0$ be equal to $J\beta^\times J^{-1}$, we find that $\epsilon \beta\beta^0=-\epsilon$, hence $\beta\beta^0=-1$.

%Now we also have $J_\beta J_\beta^{*_\beta}=1$, which yields $\beta J\beta^{-1}(\beta J)^\times \beta=\beta J\beta^{-1} J^\times \beta^2=1$. Using $\beta^{-1}=-J\beta J^{-1}$, $J^{-1}=\epsilon J$ and $J^\times=-J$ we quickly obtain $\beta^4=1$. But since $\beta^{*_\beta}=\beta$, $\beta^2$ has positive spectrum, hence $\beta^2=1$. From $\beta^0=-\beta^{-1}=-\beta$ we then get $J\beta J^{-1}=-\beta$.

Since we have $J_\beta^2=\beta J\beta J=-\epsilon\beta^2$, we quickly obtain $\beta^2=1$ as a necessary condition. 

Conversely we have already noticed that $({\cal A},*_\beta)$ is a $C^*$-algebra and   $\pi$ is a $C^*$-morphism by construction. Clearly  $D_\beta={1-i\over 2}D+{1+i\over 2}D^{*_\beta}$ is self-adjoint for $*_\beta$. Moreover, using $\chi\beta=-\beta\chi$, we have  $\chi J_\beta=-\beta \chi J=-\epsilon''(\beta J)\chi$,  and $\chi^{*_\beta}=\beta\chi^\times\beta^{-1}=-\beta\chi\beta^{-1}=\chi$.

Now we also have $J_\beta J_\beta^{*_\beta}=\beta J\beta(\beta J)^\times \beta^{-1}=\beta J\beta J^\times =-\beta^2 JJ^\times=1$  by the normalization of $\beta$.

There only remains to show that $[D_\beta,J_\beta]=0$. This calculation is completely similar to (\ref{eq16}) in \cite{part1} and left to the reader.
\end{demo}
%\bea
%[D_\beta,J_\beta]&=&\left({1-i\over 2}D+{1+i\over 2}\beta D\beta^{-1}\right)\beta J-\beta J\left({1-i\over 2}D+{1+i\over 2}\beta D\beta^{-1}\right)\cr
%&=&{1-i\over 2}D\beta J+{1+i\over 2}\beta D J-{1+i\over 2}\beta JD-{1-i\over 2}\beta J\beta D\beta^{-1}\cr
%&=&{1-i\over 2}(D\beta J+\beta\beta JD\beta)\cr
%&=&{1-i\over 2}(-DJ\beta+JD\beta)=0
%\eea 

\begin{remark}
If we did not impose $\beta$ to be imaginary in the definition of spectral spacetimes, we could say that Wick rotation to Euclidean signature is possible iff $\beta$ is both normalized and imaginary. The proof would be the same as the one of (\ref{nec12}) below.
\end{remark}

Now we try to rotate an Euclidean spectral triple $S=({\cal A},H,\pi,D,J,\chi)$ to an antilorentzian spectral spacetime. The scalar product $\bra .,.\ket$ on $H$ must be recovered as $\bra .,.\ket=(.,\omega^{-1}.)$, for a particular positive orientation form $\omega$ on the spectral spacetime. Hence $(.,.)=\bra .,\omega .\ket$. So we take a particular 1-form $\omega$ on the spectral triple as a starting point. 

\begin{propo}\label{WAL} Let $S=({\cal A},H,\pi,D,J,\chi)$ be a real even spectral triple     with KO-dimension $n$. Let   $\omega\in \Omega^1_D({\cal A},\pi)$ be a noncommutative $1$-form such that $\omega=\omega^*$ and let
$$D_\omega={1-i\over 2}D+{1+i\over 2}\omega D\omega^{-1}; (.,.)_{\omega}=\bra .,\omega.\ket, J_\omega=\omega J, \chi_{\omega}=-\chi$$

Then $S_\omega=({\cal A},H,\pi,D_\omega,J_\omega,\chi_\omega)$ is a real even    antilorentzian spacetime with KO-dimension $2-n$  on which $\omega$ is a positive 1-form if and only if 

\be
J\omega J^{-1}=-\omega,\ \omega^2=1,\mbox{ and }\omega\in \Omega^1_{D_\omega}({\cal A},\pi) \label{nec12}
\ee
%If there exists $\omega$ such that these conditions hold we say that $S$ is of WAL type, and we call $\omega$ a \emph{distinguished} 1-form on $S$. %If furthermore $\omega$ can be chosen to be exact, we say that the triple is exact. (DO WE HAVE TO REQUIRE EXACTNESS FOR BOTH COMPLEXES ?)

%\item If a spectral triple is of WAL type then for all $b\in {\cal A}$ one has $\pi(b)^0=J_\omega \pi(b)^\times J_\omega^{-1}$.

\end{propo}
  
\begin{demo} The third condition of (\ref{nec12}) is obviously necessary.  Let us first prove that the first two also are. We call $\epsilon,\epsilon''$ the KO signs of $S$. We need to have $J_\omega^2=-\epsilon$, but we have $J_\omega^2=\omega J\omega J=\epsilon\omega J\omega^* J^{-1}=\epsilon\omega\omega^0$. Hence it is necessary that $\omega\omega^0=-1$.

%From this we infer in particular that $\omega$ is invertible, and that $(.,.)_a$, which is hermitian since $\omega^*=\omega$, is non-degenerate.

Writing $\times$ for the involution corresponding to $(.,.)_\omega$, we must also have 
  $J_\omega^\times J_\omega=-1$, which yields $\omega^{-1}J_\omega^*\omega\omega J=\omega^{-1}J^*\omega^3 J=-\omega^0 \epsilon^2 J\omega^3 J^{-1}=-(\omega^0)^4$. Hence $(\omega^0)^4=1$, from which we get $\omega^4=1$.  Now $\omega^*=\omega$, thus  $\omega^2$ is a positive operator and can only be equal to $1$. From $\omega^0=-\omega^{-1}$ we obtain $J\omega J^{-1}=-\omega$, the first two conditions of (\ref{nec12}) are met.
  
Now suppose that (\ref{nec12})   holds. Then $(H,(.,.)_\omega)$ is obviously a Krein space with fundamental symmetry $\omega$. We have

\bea
D_\omega^\times&=&\left({1-i\over 2}D+{1+i\over 2}\omega D\omega^{-1}\right)^\times\cr
&=&\omega^{-1}\left({1-i\over 2}D+{1+i\over 2}\omega D\omega^{-1}\right)^*\omega\cr
&=&{1+i\over 2}\omega^{-1}D\omega+{1-i\over 2}\omega^{-2}D\omega^2
\eea

which is equal to $D_\omega$ thanks to $\omega^2=1$.

With $\beta=\omega$ we have $\beta^\times=\beta$ and $\bra .,.\ket_\beta=(.,\beta^{-1} .)_\omega=\bra . ,.\ket$ is positive definite, and for every $a\in {\cal A}$ we have $\pi(a)^{*_\beta}=\pi(a)^*=\pi(a^*)\in {\cal A}$.

It is immediate to check that $\chi_\omega^\times=-\chi_\omega$, $J_\omega\chi_\omega=-\epsilon'' \chi_\omega J_\omega$, $D_\omega\chi_\omega=-\chi_\omega D_\omega$, $J_\omega^2=-\epsilon$ and $J_\omega^\times J_\omega=-1$.
% Other obvious checks :  $\chi_\omega^\times=\omega^{-1}\chi_\omega^*\omega=\omega^{-1}\chi_\omega\omega=-\chi_\omega$, $J_\omega\chi_\omega=\omega J\chi_\omega=\epsilon''\omega\chi_\omega J=-\epsilon''\chi_\omega J_\omega$ and $D_\omega\chi_\omega=-\chi_\omega D_\omega$.
%We have $J_\omega^2=-\epsilon$ and $J_\omega^\times J_\omega=-1$ as we should by the same calculation as above. 

Finally we have:

\bea
[D_\omega,J_\omega]&=&\left({1-i\over 2}D+{1+i\over 2}\omega D\omega^{-1}\right)\omega J-\omega J\left({1-i\over 2}D+{1+i\over 2}\omega D\omega^{-1}\right)\cr
&=&{1-i\over 2}D\omega J+{1+i\over 2}\omega D J-{1+i\over 2}\omega JD-{1-i\over 2}\omega J\omega D\omega^{-1}\cr
&=&{1-i\over 2}(D\omega J+\omega\omega JD\omega)\cr
&=&{1-i\over 2}(-DJ\omega+JD\omega)=0
\eea 
%\item One has $\pi^o(b)=J\pi(b)^*J^{-1}=J\omega^2\pi(b)^*\omega^2J^{-1}=-\omega J\omega\pi(b)^*\omega(-J^{-1}\omega)=J_\omega\pi(b)^\times J_\omega^{-1}$.
\end{demo}

If there exists $\omega$ such that these conditions hold we say that \emph{$S$ is of WAL type} (for ``Wick-anti-Lorentz''), and we call $\omega$ a \emph{distinguished} 1-form on $S$. 

%Propositions \ref{ALversEUC} and \ref{WAL} taken together entail that a real even antilorentzian spectral spacetime is the Wick rotation of a real even spectral triple iff it is normalizable.
Let us observe that if a spectral triple is of WAL type then the right representation $\pi^o$ is the same for the spectral triple and its Wick rotation. %Indeed, for all $b\in {\cal A}$ one has   $\pi^o(b)=J\pi(b)^*J^{-1}=J\omega^2\pi(b)^*\omega^2J^{-1}=-\omega J\omega\pi(b)^*\omega(-J^{-1}\omega)=J_\omega\pi(b)^\times J_\omega^{-1}$.

Clearly the construction in propositions \ref{ALversEUC} and \ref{WAL} are inverse to each other. Hence the class of normalizable (and reconstructible) spectral spacetimes exactly contains the images of spectral triples under Wick rotation. In the case of manifolds, reconstructibility and normalizability are automatically satisfied. %It is not so in the noncommutative world: we will meet non-reconstructible spectral spacetimes. As for normalizability, the question is not settled yet.

\begin{remark}
We can wonder when the spectral spacetime obtained by Wick rotation of a WAL type spectral triple is stably causal. More precisely we can wonder if this happens iff the distinguished 1-form is exact. This is not obvious because exactness for $({\cal A},\pi,D)$ and $({\cal A},\pi,D_\omega)$ are not the same thing. We leave this question open.
\end{remark}

\section{Finite-dimensional case}\label{finitedim}
The considerations we have made so far, while incomplete since we have not dealt with the analytical questions, are of general validity. On the contrary, we gather in this section some results which hold only in the finite-dimensional case.

\subsection{Stably causal finite-dimensional spectral spacetimes}\label{sec51}

We begin with a simple lemma, sometimes called Jacobson's theorem, which turns out to be particularly useful.

\begin{lemma} (Jacobson's theorem) Let $A,B\in M_N(\CC)$ such that $[[A,B],A]=0$. Then $[A,B]$ is nilpotent.
\end{lemma}
\begin{demo}
Let $k\in\NN^*$. Then $\tr([A,B]^k)=\tr([A,B]^{k-1}AB)-\tr([A,B]^{k-1}BA)=\tr(A[A,B]^{k-1}B)-\tr(A[A,B]^{k-1}B)=0$ using the commutation property and the cyclicity of the trace.  
\end{demo}

Let $S=({\cal A},K,(.,.),\pi,D,J,\chi)$ be a stably causal finite-dimensional spectral spacetime, with $\beta=i[D,\pi(\delta)]$ an exact positive time-orientation form.

\begin{cor} One cannot have $\pi(\delta)^{*_\beta}=\pi(\delta)^\times$.
\end{cor}
\begin{demo}
First suppose we are in the anti-Lorentzian case. Suppose that  $\pi(\delta)^{*_\beta}=\beta\pi(\delta)^\times\beta^{-1}=\pi(\delta)^\times$. Then $\pi(\delta)^\times$ commutes with $\beta$. Thus $\pi(\delta)$ commutes with $\beta^\times=\beta$ also. But this entails that $[[D,\pi(\delta)],\pi(\delta)]=0$, which implies that $[D,\pi(\delta)]$ is nilpotent. This is absurd since   $\beta$ is Krein-positive, hence invertible.

In the Lorentzian case, we obtain that $\pi_\beta(\delta)$ commute with $\chi[D,\pi_\beta(\delta)]$. But it also commutes with $\chi$. Thus $\pi_\beta(\delta)$ commutes with $\chi^2[D,\pi_\beta(\delta)]=[D,\pi_\beta(\delta)]$, and we conclude that $[D,\pi_\beta(\delta)]$ is nilpotent. Hence $\chi[D,\pi_\beta(\delta)]$ is not invertible, a contradiction.
\end{demo}

The interest of this corollary is that since we need $\beta^\times=\beta$, an obvious guess is to look for a Krein-self-adjoint $\pi(\delta)$, and examples show that solution of this kind exist. But then $\pi(\delta)$, which plays the role of a global time function, cannot be self-adjoint for the scalar product. Hence it is not an observable in the usual quantum mechanical sense.  

\begin{cor}\label{corjac2} There are no finite-dimensional stably causal spectral spacetimes satisfying the zeroth and first-order conditions.
\end{cor}
\begin{demo}
From proposition \ref{propo12} we know that if the  zeroth and first-order conditions are satisfied then $\beta$ commutes with $\pi({\cal A})$. In particular if $\beta=i[D,\pi(\delta)]$ then $\beta$ commutes with $\pi(\delta)$ and is thus nilpotent, which is absurd.
% The first-order condition applied to $\beta=i[D,\pi(\delta)]$ shows that $[\beta,J_\beta \pi(a)J_\beta^{-1}]=0$ for all $a\in{\cal A}$. Easy manipulation then yield $\beta J\pi(a)J^{-1}=J\pi(a)J^{-1}\beta=0$, and since $\beta$ is imaginary we obtain that it commutes with $\pi({\cal A})$. By Jacobson theorem $\beta$ is nilpotent, which is absurd.
\end{demo}

\subsection{Spectral  spacetimes with algebra $\CC^2$}\label{dim2}
As a first example of finite-dimensional spectral spacetimes, let us look for all the (real, even) antilorentzian  spectral  spacetimes  with algebra ${\cal A}=\CC^2$ represented as diagonal matrices on ${\cal K}=\CC^2$.%(spinors are just scalars)%We write $\bra .,.\ket$ for the canonical scalar product on ${\cal K}$, and $j$ is matrix of the Krein product : $(.,.)=\bra .,j.\ket$. 

Since $\chi$ commutes with $\pi({\cal A})$ we have $\chi=\pm\diag(1,-1)$.

We have a non-degenerate product $(.,.)$ on ${\cal K}$ and it will be convenient to introduce its  matrix $j=\pmatrix{j_{11}&j_{12}\cr \bar j_{12}&j_{22}}$   relatively to the canonical scalar product $\bra . ,.\ket$ on ${\cal K}$, that is to say $(.,.)=\bra .,j.\ket$. As we have already noticed $\chi^\times=-\chi$ entails that $j_{11}=j_{22}=0$.  Hence we can write $j=r\pmatrix{0&e^{i\theta}\cr e^{-i\theta}&0}$ with $r>0$. We then see that $\pi(a)^\times=\pmatrix{\bar a_2&0\cr 0&\bar a_1}$, hence in this case $\pi({\cal A})$ is a $\times$-subalgebra and ${\cal A}$ is equipped with the involution:

\be
(a_1,a_2)^\times=(\bar a_2,\bar a_1)
\ee

The phase $\theta$ can be brought to $0$ by the unitary transformation $U=\diag(e^{i\theta/2},$ $e^{-i\theta/2})$ which commutes with the algebra and $\chi$, so we will assume $\theta=0$ in the following calculations.

The constraints $\chi D+D\chi=0$ and $D^\times=D$ then give $D=\pmatrix{0&b \cr c&0}$ with $b,c\in \RR$.

We easily obtain that $\Omega^1_D({\cal A},\pi)$ is generated over $\CC$ by $e_1[D,e_1]$ and $(1-e_1)[D,e_1]$ where $e_1=\diag(1,0)$. If $b=0$ or $c=0$ then $\Omega^1_D$ only contains singular matrices and there can be no time-orientation in this case. If $bc\not=0$ then $\Omega^1_D$ is the vector space of off-diagonal matrices.

Let $\beta=\pmatrix{0&\lambda\cr \mu&0}$ be a $1$-form. Then $\beta^\times=\beta$ iff $\lambda$ and $\mu$ are real. We moreover require that $(.,\beta^{-1}.)$ be positive definite, which amounts to say that $j\beta^{-1}$  is a positive definite matrix. This is so iff $\mu$ and $\lambda$ are of the same sign as $r$. Since $j\beta^{-1}$ is diagonal we evidently have $\pi(a)^{*_\beta}=\pi(a)^*$ for all $a\in{\cal A}$ (where $*$ is the adjoint for the canonical scalar product), hence $\pi(a)^{*_\beta}\in\pi({\cal A})$. The reconstructibility condition is here automatically satisfied.

We easily find that $\beta$ is exact iff ${\lambda\over b}+{\mu\over c}=0$. Since $\mu$ and $\lambda$ have the same sign, this equation has solution iff $bc<0$. When it is the case, the ``time function'' $\delta=(a_1,a_2)$ satisfying $i[D,\pi(\delta)]=i(a_2-a_1)\pmatrix{0&b\cr -c&0}=\beta$ must be such that $a_2-a_1$ is pure imaginary. 

%Hence  $ij[D,\pi(\delta)]=ir(a_2-a_1)\pmatrix{-c&0\cr 0&b}$ must be a positive matrix for some values of $a_1,a_2$, which is possible if and only if $a_2-a_1$ is pure imaginary and   $b$ and $c$ have opposite signs. 

%In summary we find the solution:
% 
%$$j=r\pmatrix{0&1\cr 1&0},\  D= \pmatrix{0& b\cr c&0},\ \chi=\pm\pmatrix{1&0\cr 0&-1}, \beta=\pmatrix{0&\lambda\cr \mu&0}$$
%
%with $r\in \RR^*$, $b,c\in\RR$ with $bc\not=0$, $r\mu>0$, $r\lambda>0$.  Moreover  the spectral spacetime is stably causal iff  $bc<0$.

%The 1-form $\beta$ is then
%$$\beta=i(a_2-a_1) \pmatrix{0&b\cr -c&0}$$
%Where $a_1,a_2$ are pure imaginary such that $i(a_2-a_1)$ has the sign of $rb$.

We know look for a real structure $J$. From $J\chi=\epsilon''\chi J$ and $J^2=\epsilon$ we obtain that

\begin{itemize}
\item If $\epsilon''=1$, then $J\psi=\diag(e^{ik_1},e^{ik_2})\bar\psi$, and $\epsilon$ can only be equal to $1$. From $JD=DJ$ we obtain that $e^{i(k_1-k_2)}=1$. However, for every Krein-self-adjoint 1-form $\beta$  we obtain $J\beta=\beta J$, hence this real structure does not define a spectral spacetime. %We can also check that $J^\times J=1$. Hence if we were to admit a real $\beta$ the  $KO$-dimension triple of signs would have been $(1,1,1)$ which  does not fit in the antilorentzian table !

%We also find that $\pi^o(a_1,a_2)=\diag(a_2,a_1)$, so $\pi^o\not=\pi$ but the condition of order $0$ is still satisfied. Also $\beta$ is real, normalizable, and the first-order condition is not fulfilled.

 %From $JD=DJ$ we obtain that $e^{i(k_1-k_2)}=e^{2i\theta}$. We then check that $J^\times J=1$. The  $KO$-dimension triple of signs is then   $(1,1,1)$ which  does not fit in the antilorentzian table !
\item If $\epsilon''=-1$ then $J\psi=\pmatrix{0&z\cr w&0}\bar\psi$ with $z\bar w=\bar z w=\epsilon$, and we find that $J^\times J=\kappa$ is fulfilled with $\kappa=\epsilon$. The axiom $JD=DJ$ ensures that $z=-w$ and $c=-b$, hence $z$ is a phase and $\epsilon$ must be equal to $-1$. Moreover we easily find that $\beta$ is imaginary $\Leftrightarrow \beta$ is exact $\Leftrightarrow \lambda=\mu$.
\end{itemize}

Hence we find a stably causal spectral spacetime with KO signs  $(-1,-1,-1)$  (this corresponds to KO dimension $2\ [8]$), such that every positive orientation form on it is exact.

We also find that $\pi^o(a)=J\pi(a)^\times J^{-1}=\pi(a)$, hence the condition of order zero is trivially satisfied.  (Clearly  the condition of order 1 is not satisfied in view of corollary \ref{corjac2}.) 

Let us summarize this solution, which we see is unique up to unitary equivalence, putting the phase back for future use:

$$j=r\pmatrix{0&e^{i\theta}\cr e^{-i\theta}&0},\  D=b \pmatrix{0& e^{i\theta}\cr -e^{-i\theta}&0}, \chi=\pmatrix{-1&0 \cr 0&1}, J\psi= \pmatrix{0&-1\cr 1&0}\bar\psi$$

where we have chosen the overall sign of $\chi$ and  phase of $J$, which play no role, in order to ease comparison with other examples to be defined below.%, and $\lambda$ and $\mu$ are of the sign of $r$.

%If $\beta$ is a given positive time orientation, then the condition $\beta\omega^{-1}>0$ on imaginary 1-forms $\omega$ selects those which are also positive time orientations, which corresponds, dually, to future-pointing vector fields.

We can understand the role of the time-orientation in the following way: if $\beta=[D,f]$ where $f$ is a real function on the set $\{1;2\}$, then $[D,f]$ is a positive time-orientation iff $rb(f_2-f_1)>0$. Hence the sign of $rb$ determines an  ordering on $\{1;2\}$, and $[D,f]$ is a positive time-orientation iff $f$ is increasing for this ordering.

% the real function $g$ defines through $[D,g]$ another positive time-orientation iff $f$ and $g$ have the same sense of variation. Hence a positive time-orientation induces a particular   ordering on $\{1;2\}$. Changing the sign of either $r$ or $b$ reverses this ordering.

%, since $i(a_2-a_1)$ is real and $J$ anticommutes with $\pmatrix{0&b\cr b&0}$. Hence the condition of order 1 is not satisfied.

We remark that Euclidean even real spectral triple over the algebra ${\cal A}=\CC^2$ and with Hilbert space $\CC^2$  come in two breeds, up to unitary equivalence:

\begin{itemize}
\item One of KO dimension $0$, with Dirac operator $D=b\pmatrix{0&1\cr  1&0}$, $\chi=\pmatrix{1&0\cr 0&-1}$, $J=e^{ik}\bar \psi$, with $b>0$, $k\in\RR$. We call it $S^0$.    %one of KO dimension $0$, with Dirac operator $D=|b|\pmatrix{0&e^{i\theta}\cr  e^{-i\theta}&0}$, $\gamma=\pmatrix{1&0\cr 0&-1}$, $J=\pmatrix{e^{ik_1}&0\cr 0&e^{ik_2}}\circ$ c.c., with $e^{i(k_1-k_2}=e^{2i\theta}$,  
\item One of KO dimension $6$, with the same $D$ and $\chi$ but with $J=e^{i\phi}\pmatrix{0&1\cr 1&0}\circ$ c.c. We call it $S^6$.
\end{itemize}

%It is easily seen that the first is the Wick rotation of the antilorentzian spectral spacetime we found above. 

It turns out that $S^0$ is of WAL type while $S^6$ is not. Indeed, one easily verifies that a 1-form $\omega$ on $S^0$ satisfies $\omega^*=\omega$, $\omega^2=1$ and $J\omega=-\omega J$ iff it is of the form $\omega=\pm\pmatrix{0&i\cr -i&0}$. Choosing one of these two we find that $D_\omega={1-i\over 2}D+{1+i\over 2}\omega D\omega=-iD$. Hence $\omega\in \Omega^1_{D_\omega}$, and $S^0$ is of WAL type. The Wick rotated spectral spacetime $(S^0)_\omega$ corresponds to the solution we found above with $\theta=-\pi/2$.

%The ``global time functions'' $\delta\in{\cal A}$ such that $i[D,\pi(\delta)]=\omega=\pmatrix{0&i\cr -i&0}$ are of the form $diag(a_1,a_2)$ with $a_2-a_1=i$. 

%The element $a\in{\cal A}$ (the ``global time function'') which permits to Wick rotate it to antilorentzian signature is any one of the form $a=\diag(x,y)$, with $x-y\in\RR$ and $|x-y||b|=1$. The Wick rotated Dirac operator is $D_\omega=-iD+$ diagonal term. Hence we see that the element $\delta$ such that $[D_\omega,\pi(\delta)]=[D,\pi(a)]$ must be of the form $\delta=ia+$ a term commuting with $D$. Of course we see that we recover the antilorentzian spectral triple with global time function of KO dimension $s-t=2 \ [8]$ in this way.

%We see an interesting phenomenon: the global time function $\delta$ for the antilorentzian spectral triple cannot be self-adjoint ! This might seem strange at first, but notice that instead of self-adjointness we can have Krein-adjointness : $\delta^\times=\delta$ if we choose $ \delta=diag(-i/2,i/2)$. (Note that in this case $\delta$ is real in the sense that it commutes with $J$.) The consideration of General Relativistic Quantum Mechanics will show that it is exactly as it should be. Time observables are not observable in the usual QM sense (see below).

The reason why  $S^6$ is not of WAL type is that a self-adjoint 1-form $\omega$ necessarilly commutes with $J$.% If we try to Wick rotate $S^6$ nonetheless, using a real form instead of a imaginary one, we obtain the first case which has been eliminated above, with KO signs $(1,1,1)$ not fitting in the table, and which thus be described as a ``bad Wick rotation''.

\subsection{The canonical triple of a positively weighted graph and its Wick rotation}\label{sec64}

The   finite analogue of a Riemannian manifold is a finite undirected simple graph with positive weights on the edges. More precisely we consider a set of vertices $V$, with  $|V|=n$, a set of edges $E$ and a weight function $\delta : E\rightarrow \RR_+^*$. An element $e\in E$ is a pair $\{i,j\}\subset V$ with $i\not=j$. This naturally induces a geodesic distance $d$ on the vertices in much the same way as in Riemannian manifolds: $d(i,j)$ is the infimum of the lengths of paths joining $i$ and $j$. We allow disconnected graphs, so that this distance can reach infinity.

We look for a representation of the $C^*$-algebra ${\cal A}={\cal C}(V)\simeq \CC^n$ on a Hilbert space and a Dirac operator $D$ such that  $d$ is recovered through Connes' distance formula. That is, we want:

\be
d(i,j)=\sup\{|a(i)-a(j)|\ |\  a\in{\cal A}, \|[D,\pi(a)]\|\le 1\}
\ee

A first guess would be to set $H=L^2(V)$, let ${\cal A}$ act diagonally and set $D_{ij}={1\over \delta_{ij}}$, where for convenience we have defined $\delta_{ij}:=\delta(\{i;j\})$. Thus $D$ would just be a generalized adjacency matrix of the  graph. Unfortunately this solution does not work, except in very special circumstances, for instance when $n=2$. The correct solution, which we will describe now, is well-known (see  \cite{ikm} or \cite{vs} p 19). We depart a little from the usual presentation by allowing the graph to be incomplete and the weights not be equal to the distance (that is, $\delta_{ij}$ is not necessary equal to $d(i,j)$). % We make a small change by considering a representation twice as much degenerate as usual, which will allow us to consider general directed graphs later on.

%Let $E$ be the set of  couples $(i,j)\in V\times V$ such that $0<\delta_{ij}<\infty$. The source and target functions $s$ and $t$ are just $s(i,j)=i$, $t(i,j)=j$. Then $(V,E,s,t)$ is a simple directed graph. 

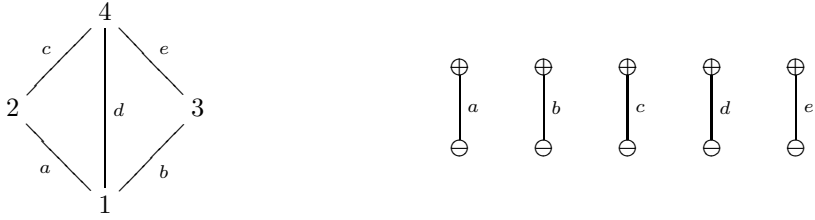
\begin{figure}[hbtp]
%\begin{center}
\parbox{6cm}{$\xymatrix{ & 4 & \cr
2\ar@{-}[ur]^c & & 3\ar@{-}[ul]_e\cr
 & 1\ar@{-}[ul]^a \ar@{-}[uu]_d \ar@{-}[ur]_b &
}$}\parbox{7cm}{$\entrymodifiers={[o][F-]}\xymatrix{+&+&+&+&+\cr
-\ar@{-}[u]_a&-\ar@{-}[u]_b&-\ar@{-}[u]_c&-\ar@{-}[u]_d&-\ar@{-}[u]_e}$}
%\end{center}
\caption{A undirected weighted graph $G=(V,E,\delta)$ on the left and its split form $\tilde E$ on the right (wieghts are not represented).}\label{splitgraph}
\end{figure}

We set $\tilde E=E\times \{-,+\}$, where $\{-,+\}$ is a two-element set. Elements of $\tilde E$ are to be interpreted as the extremities of the edges of the graph when it is splitted into a disjoint union of its $1$-edge subgraphs (see figure \ref{splitgraph}). This is why we call $\tilde E$ \emph{the split graph} of $G$.

We then take $H=L^2(\tilde E)\simeq \CC^E\otimes\CC^2$, with the canonical scalar product.  In order to define our representation we will choose sources and targets for each edge, that is we will give an arbitrary orientation to our graph which then becomes a simple directed graph. The source and target functions are called $s,t : E\rightarrow V$, and it will be convenient to use the notation $e^+:=t(e)$ and $e^-=s(e)$.  We define the representation $\pi$ by

\be
 (\pi(a)F)(e,\pm)=a(e^\pm)F(e,\pm)\nonumber
\ee

where $a\in {\cal A}$, $F\in H$ and $e$ is an edge.

%The Dirac operator is generally defined simply by $(DF)(e,\pm)={1\over \delta_e}F(e,\mp)$, where we have introduced the notation $\delta_e=\delta(s(e),t(e))$. However we want to be as general as possible, 

The Dirac operator is  defined simply by:

\be
(DF)(e,\pm)={e^{\pm i\theta_e}\over \delta_e}F(e,\mp)\nonumber
\ee

The phases $\theta_e$ can be removed by a unitary transformation but we keep them in order to make the comparison easier with the antilorentzian situation.

Using the isomorphism $\CC^E\otimes \CC^2\simeq \bigoplus_{e\in E}\CC^2$ and the basis $(-,+)$ of $\CC^2$ we can write down our definitions matricially as

\be
\pi(a)=\bigoplus_{e\in E}\pmatrix{a(e^-)&0\cr 0&a(e^+)},\ D=\bigoplus_{e\in E}{1\over\delta_e}\pmatrix{0&e^{i\theta_e}\cr e^{-i\theta_e}&0}
\ee

The commutator $[D,\pi(a)]$ comes out as

\be
[D,\pi(a)]=\bigoplus_{e\in E}{a(e^+)-a(e^-)\over \delta_e}\pmatrix{0&e^{i\theta_e}\cr -e^{-i\theta_e}&0}\label{comutdpia}
\ee

so that $\|[D,\pi(a)]\|=\sup_{e\in E}{|a(e^+)-a(e^-)|\over \delta_e}$. If $\|[D,\pi(a)]\|\le 1$ and $(i_0,\ldots,i_k)$ is a geodesic path in the (undirected) graph joining $i_0$ and $i_k$, then we see that 

\bea
|a(i_k)-a(i_0)|&\le&\sum_{r=0}^{k-1}|a_{i_{r+1}}-a_{i_r}|\cr
&\le& \sum_{r=0}^{k-1}\delta_{r,r+1}\mbox{ by }(\ref{comutdpia})\cr
&=&d(i_0,i_k)
\eea
since $(i_0,\ldots,i_k)$ is a geodesic path. Moreover if we take $a$ to be the function $a(j):=d(i_0,j)$ then we have $a(i_k)-a(i_0)=d(i_0,i_k)$, so that we just need to prove that $\|[D,\pi(a)]\|\le 1$ for this particular function. This amounts to prove that for each edge $e=\{j_1,j_2\}$ one has $|d(i_0,j_1)-d(i_0,j_2)|\le \delta_{j_1,j_2}$. But this is true because $|d(i_0,j_1)-d(i_0,j_2)|\le d(j_1,j_2)$ by the triangle inequality and $d(j_1,j_2)\le \delta_{j_1,j_2}$.

An obvious grading which commutes with the representation of the algebra and anti-commutes with $D$ is given by

\be
(\chi F)(e,\pm)=\pm F(e,\pm)\label{defgamma}
\ee

This is in fact (up to an overall sign) the only sensible possibility for $\chi$. Indeed, from $\chi D^2=D^2\chi$ we see that, provided the weights $\delta_e$ are all different, $\chi$ must have the form $\bigoplus_e \chi_e$, with $\chi_e$ commuting with the restriction of $\pi({\cal A})$ to the edge $e$. Thus (\ref{defgamma}) is the only definition for $\chi$ for which the signs do not depend on the edges and  is thus consistent with every simple directed graph we can start with.

Hence we see that if there is only one edge we recover   the most general Euclidean spectral triple with algebra $\CC^2$. We know then that there are two possible real structures that we could put on our spectral triple, giving respective KO-dimension $0$ and $6$ (see section \ref{dim2}). However we have already seen that the second one cannot be Wick rotated to antilorentzian signature. Since this is what we have in mind, we define

\be
JF(e,\pm)=e^{\mp i\theta_e}\bar F(e,\pm)
\ee

\be
J=\bigoplus_{e\in E}\pmatrix{e^{i\theta_e}&0\cr 0&e^{-i\theta_e}}\circ c.c
\ee

Hence we obtain an even real spectral triple $S^G=({\cal A},\pi,H,D,J,\chi)$ of KO-dimension $0$ such that the distance function on the original metric space $V=\widehat{\cal A}$ is recovered by  Connes' distance formula. We will call it the \emph{canonical triple} of the weighted graph $G$ (even though the arbitrary phases and orientation make it non-unique, all these choices define the same unitary equivalence class of spectral triples).

We will now show that it is of WAL type.

\begin{propo}  The canonical triple over a weighted undirected graph $G$ is of WAL type.
%Let $\omega\in\Omega^1_D({\cal A},\pi)$ be a self-adjoint imaginary and normalized  1-form. Then $\omega$ induces an orientation on the edges of $G$. %This orientation is acyclic if $\omega$ is exact.
%\item Let $S_\omega^G$ be the antilorentzian spectral spacetime obtained by Wick rotation of $S$ through $\omega$. Then $S_\omega^G$ is stably causal iff the orientation induced by $\omega$ is acyclic.
\end{propo}
\begin{demo}
%For this let us observe thanks to (\ref{comutdpia}) that a general 1-form $\sum_i \pi(a_i)[D,\pi(b_i)]$ can be written as

First let us calculate $\Omega^1_D(\pi,{\cal A})$. Let $e=\{i,j\}$ be an edge, with $s(e)=i$ and $t(e)=j$. Let $\chi_i$ be the function which takes value $1$ on $i$ and $0$ on the other vertices. Then $\pi(\chi_j)[D,\pi(\chi_i)]=\bigoplus_{e'}A_{e'}$ where $A_{e'}$ is a matrix which is zero except for $e=e'$ and $A_e=\pmatrix{0&?\cr 0&0}$ with $?$ a nonzero coefficient. Similarly $\pi(\chi_i)[D,\pi(\chi_j)]$ contains only a nonzero summand for the edge $e$ which is of the form $\pmatrix{0&0\cr ?&0}$. Since $\Omega^1_D(\pi,{\cal A})$ is an ${\cal A}$-bimodule, we easily obtain that the general expression of $\omega$ is $\omega=\bigoplus_{e\in E}\pmatrix{0&\omega_e\cr \omega_e'&0}$. We can then  check that $\omega$ is imaginary iff $\omega_e e^{i\theta_e}$ and $\omega_e' e^{-i\theta}$ are imaginary for all $e$. Since we want in addition $\omega$ to be   self-adjoint it must be of the form $\omega=\bigoplus_{e\in E}\pmatrix{0&ix_e e^{i\theta_e}\cr -ix_ee^{-i\theta_e}&0}$ with $x_e\in\RR$ and finally since it is normalized we have

\be
\omega=\bigoplus_{e\in E}\sigma_e\pmatrix{0&ie^{i\theta_e}\cr -ie^{-i\theta_e}&0}\label{eq19}
\ee

where $\sigma_e$ is a sign. We then find that 

\be
D_\omega=\bigoplus_{e\in E} {i\over\delta_e}\pmatrix{0&e^{i\theta_e}\cr e^{-i\theta_e}&0}=-iD
\ee

Thus $\Omega_{D_\omega}^1({\cal A},\pi)=\Omega_D^1({\cal A},\pi)$ and we see that all the conditions of proposition \ref{WAL} are satisfied, hence the canonical triple is of WAL type.
\end{demo}

The  signs $\sigma_e$ define in an obvious way an orientation on the graph $G$.  For ease of notation we  redefine the source and target functions which were arbitrarily chosen to define $\pi$ in order to make the $\sigma_e$ be all equal to $+1$. Hence we suppose that the fiducial orientation   defined by $s,t$ and the one defined by $\omega$ coincide.

Then we find that

\be
J_\omega=\bigoplus_{e\in E}\pmatrix{0&-1\cr 1&0}\circ c.c\label{eq21}
\ee

and 
\be
(F,G)_\omega=\sum_{e\in E}\left((F(e,-),ie^{i\theta_e}G(e,+))+(F(e,-),-ie^{-i\theta_e}G(e,+))\right)\label{eq22}
\ee

The spectral spacetime $S^G_\omega=({\cal A},H,(.,.)_\omega,\pi, D_\omega, J_\omega,\chi)$ can clearly be defined directly by (\ref{eq19})-(\ref{eq22}), without reference to the spectral triple $S^G$,  as soon as the orientation on $G$ defined by $\omega$ is given. Hence we can associate an  antilorentzian spectral spacetime to any oriented wieghted graph, which we call the \emph{canonical antilorentzian spectral spacetime} over this structure. Let us now consider the question of its stable causality.
 
Let  $\beta=\bigoplus_e\pmatrix{0&\beta_e\cr  \beta_e'&0}$ be a  1-form. It  anticommutes with $J_\omega$ iff $\beta_e=\bar \beta_e'$ and it is Krein self-adjoint iff $\beta_e e^{-i\theta_e}$ and $\beta_e'e^{i\theta_e}$ are pure imaginary.  Hence to have both properties we must have $\beta_e=ix_e e^{i\theta_e}$ with $x_e\in\RR$ and $\beta_e'=-ix_e e^{-i\theta_e}$ as above. Then $\beta$ will be a positive time-orientation for $S_\omega^G$ iff $\beta\omega=\bigoplus_e x_e I_2$ is a positive operator, which means that all $x_e$ are positive.

\begin{figure}[hbtp]
%\begin{center}
\hspace{2cm}\parbox{5cm}{$\xymatrix{ & 4\ar[dr]^e& \cr
2\ar[ur]^c & & 3\ar[dl]_b\cr
 & 1\ar[ul]^a \ar[uu]_d &
}$}\parbox{6cm}{$\xymatrix{ & 4 & \cr
2\ar[ur]^c & & 3\ar[ul]_e\cr
 & 1\ar[ul]^a \ar[uu]_d \ar[ur]_b &
}$ }
%\end{center}
\caption{On the left $\omega$ is chosen so as to induce directed cycles, such as $(1,4,3,1)$. The resulting triple $S_\omega^G$ is not stably causal. On the right $\omega'$ induces an acyclic orientation. The triple $S^G_{\omega'}$ is stably causal, a strictly increasing function on the vertices of the graph being the identity. However $\omega'$ is not itself exact: the integral of $\omega'$ on the oriented cycle $(1,4,3,1)$ being$+1-1-1\not=0$.}\label{fig2}
\end{figure}
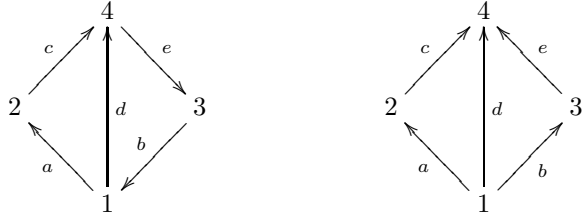

Let ${\cal P}=(v_1,\ldots,v_k)$ be a path in the graph, that is a finite sequence of vertices such that $v_{i+1}$ is adjacent to $v_i$. Then define 

\be
\int_{\cal P}\beta=\sum_{i=1}^{k-1}\epsilon_{i,i+1} x_{v_i,v_{i+1}}\delta_{\{v_i;v_{i+1}\}}\label{blue}
\ee
 
where $\epsilon_{i,i+1}=1$ if $v_i$ is the source of the edge $\{v_i;v_{i+1}\}$ and $\epsilon_{i,i+1}=-1$ otherwise. We say that ${\cal P}$ is a cycle if $v_1=v_k$. Note that the orientation of ${\cal P}$ defined by the signs $\epsilon_{i,i+1}$ need not correspond to the one induced by $G$ (see figure \ref{fig2}).

\begin{remark}
Definition (\ref{blue}) can be rephrased in a way which may be more natural in this context. Consider ${\cal P}$ as a subgraph of the undirected graph $G$, build the canonical triple over it, and consider the orientation 1-form $P=\bigoplus_{e\in {\cal P}}\epsilon_e \pmatrix{0&ie^{i\theta_e}\cr -ie^{-i\theta_e}&0}$. Then
$$\int_{\cal P}\beta=\tr(\beta P|D_{\cal P}|^{-1})=(\beta,P)$$

where $(\omega_1,\omega_2)=\tr_\omega(\omega_1^\times\omega_2|D|^{-d})$ is the obvious generalization to the semi-Riemannian setting of the usual scalar product on $1$-forms defined on a spectral triple (\cite{redbook}, chap. VI). Note that one has to put the dimension $d=1$ by hand in this formula.
\end{remark}

We can now state a  discrete analogue of Morera's theorem.

\begin{lemma} The form $\beta$ is exact iff $\int_{\cal P}\beta=0$ for any cycle ${\cal P}$.
\end{lemma}
\begin{demo} Let $\beta=i[D,f]$ be an exact 1-form. We see that from (\ref{comutdpia})  that $x_e={f(e^+)-f(e^-)\over \delta_e}$, thus 

\bea
\int_{\cal P}i[D,f]\omega&=&\sum_{i=1}^{k-1}\epsilon_{i,i+1} (f(t(\{v_i;v_{i+1}\}))-f(s(\{v_i;v_{i+1}\})))\cr
&=&\sum_{i=1}^{k-1}  (f(v_{i+1})-f(v_i))\cr
&=&f(v_k)-f(v_1)=0,\mbox{ since }{\cal P}\mbox{ is a cycle}
\eea

Conversely, if $\beta$ satisfies the stated property, let us take a vertex $v_0$ and define 

\be
f(v)=\int_{v_0\rightarrow v}\beta\omega
\ee

where $v_0\rightarrow v$ is any path from $v_0$ to $v$. This definition clearly does not depend on the chosen path. If $e$ is an edge we see that $f(e^+)-f(e^-)=x_e\delta_e$, so that $i[D,f]=\beta$.
\end{demo}

Now let $\beta$ be a positive time-orientation. Since all $x_{e}$ are positive it is clear from ``Morera's theorem'' that $\beta$ cannot be exact if the graph orientation  induces any directed cycle. In that case $S_\omega^G$ is not stably causal.

Now it is a well-known fact that if $G$ is directed acyclic then there exists a topological ordering on its vertices, that is to say a strictly increasing function for the partial order defined by the graph (see for instance \cite{clrs}, pp 549--552). Let $f$ be such a function. Then $[D_\omega,\pi(f)]$ is an exact positive time-orientation for $S_\omega^G$.

Hence we proven that:

\begin{propo} The canonical spectral spacetime on a directed weighted graph is stably causal iff the graph is acyclic.
\end{propo}

%Thanks to the previous theorem we can attach a spectral spacetime to a given simple directed weighted graph $G$: we consider the canonical spectral triple on the underlying undirected weighted graph, we take $\omega$ the (unique) normalized imaginary and self-adjoint 1-form which defines the correct orientation, and Wick rotate to antilorentzian signature thanks to $\omega$. The resulting spectral spacetime can of course be defined directly, without reference to any Euclidean spectral triple. We will call it \emph{the canonical (antilorentzian) spectral spacetime over $G$}.

Let us briefly conclude this section. The canonical spectral spacetime over a graph is a very natural object for several reasons: 1) since the split graph and the original graph have the same edges, a 1-form can be pictured as a field of covectors defined on them, 2) when we Wick rotate it using a distinguished 1-form we recover a spectral triple with a natural distance function, and  3) the stable causality condition is just the existence of a partial order on the vertices compatible with the orientation of the graph.  However, we see immediately that no time-orientation form commutes with $\pi({\cal A})$, and that for $a\in{\cal A}$ we have $\pi(a)^\times=\omega\pi(a)^*\omega=\bigoplus_e \pmatrix{\overline{a(e^+)}&0\cr 0&\overline{a(e^-)}}$ and this does not belong to $\pi({\cal A})$ in general. Hence ${\cal A}$ is not a $\times$-subalgebra. For this reason, this example is not covered by the approaches to Lorentzian noncommutative geometry reviewed in section \ref{sec44}.

\subsection{A noncommutative example: the split Dirac structure}\label{sec54}
\subsubsection{Definition}
We would now like to build a genuinely noncommutative spectral spacetime over a positively weighted oriented graph $G$, which would generalize the canonical spectral spacetime over $G$. A natural idea would be to build a discrete Clifford bundle and take for ${\cal A}$ the algebra of sections of this bundle. We list below the first ingredients at our disposal  to ease future reference.

\begin{itemize}
\item $G=(V,E,s,t,\delta)$ is a positively weighted oriented and connected graph, $V$ its set of vertices, $E$ its set of edges, $s,t : E\rightarrow V$ the source and target functions. We still write $e^-=s(e)$ and $e^+=t(e)$.
\item For each $v\in V$, $X_v$ is a $n$-dimensional real vector space, and $g_v$ is an antilorentzian bilinear form on it.  
\item ${\cal A}$ is the  algebra of   sections $a : v\mapsto a_v\in \CC l(X_v,g_v)$ of the discrete Clifford bundle $\bigsqcup_{v\in V}\CC l(X_v,g_v)$. We write $c_v$ for the canonical real structure which permits to recover $Cl(X_v,g_v)$ in $\CC l(X_v,g_v)$. % We define the corresponding adjunction  $\times$  by $a^\times(v)=a(v)^\times=c_v(a(v)^T)$.
 
\item We let $S_{v}$ be an irreducible representation space of $\CC l(X_{v},g_v)$. The representation map $\rho_v$ will be often understood. We fix a $c_{v}$-compatible Krein product $(.,.)_v$  on $S_{v}$ which we often write simply $(.,.)$     to avoid notation clutter. The future cone in $X_v$ is defined to be the set of timelike vectors which turn $(.,.)_v$ into a positive definite product. %We set $K=\bigoplus_{e^\pm\in\tilde E}S_{e^\pm}$, which we see as the space of section of a discrete spinor bundle. Hence we write its elements as $F : e^\pm\mapsto F(e,\pm)\in S_{e^\pm}$. Most of the time we write the action directly by multiplication, otherwise we call $\rho_e^\pm$ the representation. 
  
% \be
% (F,G)=\sum_{e\in E}r_e\big((F(e,+),e^{-i\theta_e}h_e^+G(e,-))+(F(e,-),e^{i\theta_e}h_e^-G(e,+)\big)
% \ee

%We have to require that $r_e\not=0$ for this product to be non-degenerate. 

%\end{itemize}

%We already know that in order to build a meaningful\footnote{That is, which corresponds in some way to the geodesic distance on the graph.} Dirac operator we will have to introduce the split graph $\tilde E=E\times \{-;+\}$. 
%But let us pretend for a moment that we do not know, or do not care, about this problem. The most natural thing to do would be to mimic the construction of the canonical triple over a manifold by introducing the spinor bundle ${\cal S}=\bigcup_{v\in V}S_v$ and its space of sections  $\Gamma {\cal S}=\bigoplus_{v\in V}S_v$, equipped with the Krein product $(\Psi,\Phi):=\sum_{v\in V}(\Psi(v),\Phi(v))$. But this could not work, since the grading $\gamma$ cannot at the same time commute with ${\cal A}$ and satisfy $\gamma^\times=-\gamma$ for the (diagonal) Krein product $(.,.)$. A possible way out would be to replace ${\cal A}$ by the even subalgebra ${\cal A}^0$, but a more natural solution is provided by the split graph.

%Let us now list the remaining ingredients of the construction. 

%\begin{itemize}
\item $\tilde E=E\times\{-;+\}$ is the ``split graph'' of $G$.
\item We let  $K$ be the space of sections  $F : (e,\pm)\mapsto F(e,\pm)\in S_{e^\pm}$ of the discrete spinor bundle $\bigcup_{(e,\pm)\in\tilde E}S_{e^\pm}$ over $\tilde E$.   %We set $K=\bigoplus_{e^\pm\in\tilde E}S_{e^\pm}$, which we see as the space of section of a discrete spinor bundle. Hence we write its elements as $F : e^\pm\mapsto F(e,\pm)\in S_{e^\pm}$. Most of the time we write the action directly by multiplication, otherwise we call $\rho_e^\pm$ the representation. 
\item The representation of ${\cal A}$ on $K$ is defined by:

\be
(\pi(a)F)(e,\pm)=a(e^\pm)\cdot F(e,\pm)
\ee

\item We fix discrete ``parallel transport'' operators $h_e^+ : S_{e^-}\rightarrow S_{e^+}$, and write  $h_e^-=(h_e^+)^{-1}$. The set of all these operators will be called the \emph{discrete connection} in what follows. 
\item The Krein product on $K$ is defined to be:

\be
(F,G)=\sum_{e\in E} \big((F(e,+),h_e^+G(e,-))+(F(e,-),h_e^-G(e,+)\big)\label{defkrein}
\ee

% \be
% (F,G)=\sum_{e\in E}r_e\big((F(e,+),e^{-i\theta_e}h_e^+G(e,-))+(F(e,-),e^{i\theta_e}h_e^-G(e,+)\big)
% \ee

%We have to require that $r_e\not=0$ for this product to be non-degenerate. 
\item The Dirac operator $D$ is defined by:

\be
(D F)(e,\pm)={\mp  1\over \delta_e}\gamma_e^\pm h_e^\pm F(e,\mp)\label{Dirac49}
\ee
 
% \be
% (D F)(e,\pm)={\mp e^{\mp i\theta_e}\over \Delta t_e}\gamma_e^\pm h_e^\pm F(e,\mp)
% \ee
 
where $\gamma_e^\pm$ is   an element of $End(S_{e^\pm})$  and $\delta_e\in]0;+\infty[$.  Of course $\delta_e$ could be absorbed in the definition of $\gamma_e^\pm$ but in this way we will see more clearly the connection with the discretization of the usual Dirac operator.
\item The chirality $\chi$ is defined by:
\be
(\chi F)(e,\pm)=\chi_{e^\pm} F(e,\pm)
\ee
where $\chi_{e^\pm}$ is the chirality operator of $\CC l(X_{e^\pm})$. 
\item The real structure $J$ is:

\be
(JF)(e,\pm)= h_e^\pm J_{e^\mp} F(e,\mp)
\ee
where $J_{e^\pm}$ is an antilinear operator on $S_{e^\pm}$ satisfying $(J_{e^\pm})^2=\epsilon$, $J_{e^\pm} (J_{e^\pm})^\times=\kappa=-1$ and $c_{e^\pm}(a)=-J_{e^\pm} a J_{e^\pm}$ for any $a\in \CC l(X_{e^\pm},g_{e^\pm})$.
\end{itemize}

Some remarks about the previous definitions are in order. 

\begin{enumerate}
\item We immediately see that $\chi$ does not commute with the algebra ${\cal A}$, so we must restrict to the subalgebra ${\cal A}^0$ of even elements.
\item We note that the above definition reproduces, as is intended, the canonical spectral spacetime over a weighted directed graph when the dimension of $X_v$ is $0$,  with all the phases $\theta_e$ being equal to $-\pi/2$, $\chi_{e^\pm}=\pm 1$, $J_{e^\pm}=\pm c.c.$,  $h_e^\pm=1$, and $\gamma_e^\pm=1$. 
\item The $n=0$ case is the only one for which $S_{e^\pm}$ is an irreducible ${\cal A}^0$-module. For $n\ge 2$ we will denote $S_{e^\pm}^{L/R}$ for the $-1/+1$-eigenspaces of $\chi_{e^\pm}$, respectively. Each of them is an irreducible ${\cal A}^0$-module.
\item The algebra ${\cal A}^0$ is noncommutative for $n\ge 4$.
\end{enumerate}

We call the structure $S=({\cal A}^0,K,\pi,D,J,\chi)$ defined above \emph{ the split Dirac structure over}   $G$. Now we ask under what conditions a split Dirac structure is a spectral antilorentzian spacetime and what are its KO signs  $\tilde\epsilon,\tilde\epsilon''$. We will see   conditions on the discrete connection   emerge which are very similar to the ones discussed in the continuous case in \cite{part1}, section \ref{sec33}. We thus need to set up some terminology about the discrete connection.

\subsection{The discrete connection and holonomy group of a discrete spinor bundle}

The notations being the same as in the previous subsection, we begin with the definition of the discrete analogues of the properties reviewed in \cite{part1}, section \ref{sec33}.

\begin{definition} 
\begin{enumerate}
%\item If for all edge $e$, the parallel transport operators $h_e^\pm$ intertwine the action of $\CC l(X_{e^-})$ and $\CC l(X_{e^+})$, the discrete connection is said to be \emph{spin-c}.
\item If the parallel transport operators are Krein-unitary transformations, that is: $(h_e^+s_-,h_e^+s_-)=(s_-,s_-)$, then we say that the discrete connection is \emph{metric}.
%\item They intertwine the action of vectors at $e^-$ with the action of vectors at $e^+$, that is: for every vector $x^-\in X_{e^-}$, there exists a vector $x^+\in X_{e^+}$ such that  $\rho_{e^-}(x^-)=h_e^-\rho_{e^+}(x^+)h_e^+$.
%\item The map $x^-\mapsto x^+$ (which is necessarily an isometry since $(x^-)^2=(x^+)^2$) is asked to send the future cone at $e^-$ to the future cone at $e^+$.
\item If for all $e$, $h_e^\pm J_{e^\pm}=J_{e^\mp} h_e^\pm$ we say that the discrete connection is \emph{spin preserving}.
\item If for all $e$, $\chi_{e^\pm} h_e^\pm=h_e^\pm \chi_{e^\mp}$, we say that the discrete connection  \emph{preserves orientation}.
\end{enumerate}
\end{definition}

Note that on a space and time oriented manifold, spin connections automatically preserve space and time orientations by continuity. This is no longer true in the discrete case, hence these properties must be separately imposed. Moreover, since we have no notion of Levi-Civita transport from $X_{e^-}$ to $X_{e^+}$, we cannot directly translate the definition of Clifford connections from the continuous case to the discrete case.  Instead we will say the following:

\begin{definition}\label{defcliffcon}
The discrete connection is a \emph{Clifford connection} iff for every edge $e$ there exists a linear map $\Lambda :  X_{e^-} \rightarrow X_{e^+}$  such that the following diagram commutes for all $v\in X_{e^-}$:

$$\xymatrix{S_{e^-}\ar[r]^v\ar[d]_{h_e^+}& S_{e^-}\ar[d]^{h_e^+}\cr S_{e^+}\ar[r]^{\Lambda(v)}& S_{e^+}}$$
\end{definition}

Of course the similar property with $S_e^+$ and $S_e^-$ exchanged is then automatically satisfied with $\Lambda^{-1}$ instead of $\Lambda$. The discrete connection will be said to be \emph{a spin connection} if it is  Clifford, metric, spin and orientation preserving.

Let us see how a discrete spin connection permits to define the holonomy group of the graph. Let ${\cal P}=(x_1,\ldots,x_k)$ be a path from $x_1$ to $x_k$. If $(x_i,x_{i+1})$ is an edge $e_i$ of $G$ we set $\sigma_i=+$, if $(x_{i+1},x_i)$ is an edge we set $\sigma_i=-$. We then define the parallel transport operator $h_{\cal P} : S_{x_1}\rightarrow S_{x_k}$ to be 

$$h_{\cal P}:=\prod_{i=1}^{k-1}h_{(x_i,x_{i+1})}^{\sigma_i}$$

We also define $H_{\cal P}: \CC l(X_{x_1})\rightarrow \CC l(X_{x_k})$ to be the map $a\mapsto h_{\cal P}a h_{\cal P}^{-1}$. If ${\cal P}$ is a closed loop we write ${\rm hol}_{\cal P}:=h_{\cal P}$, and call it the \emph{spinor-holonomy} operator around ${\cal P}$. We also write ${\rm Hol}_{\cal P}:=Ad_{{\rm hol}_{\cal P}}$ and call it the \emph{holonomy} operator around ${\cal P}$. We can see now that definition \ref{defcliffcon} is very natural since it ensures that the spinor-holonomy around any loop  belongs to the Clifford group. Moreover, note that  $\Lambda(v)^2=h_e^+ v^2 h_e^-=v^2$, hence $\Lambda$ is an isometry. If in addition the discrete connection is metric, then $(\psi,\Lambda(v)\psi)=(h_e^-\psi,v h_e^-\psi)$, thus $\Lambda(v)$ will be future-directed timelike if, and only if, $v$ is, thus $\Lambda$ will preserve time-orientation. Finally if the discrete connection preserves orientation, then for any positively oriented pseudo-orthonormal basis $(e_1,\ldots,e_n)$ , $h_e^+e_1\ldots e_n h_e^-=\Lambda(e_1)\ldots\Lambda(e_n)$ will be positively oriented. Therefore $\Lambda$ will preserve the total orientation. It then makes sense to \emph{define} the Levi-Civita parallel transport along $e$ to be this $\Lambda$.  Thus we see that the spinor-holonomy around a loop of a spin connection is an element of the neutral component of the spin group, and that the corresponding holonomy operator is   the canonical extension to the Clifford algebra of a proper orthochronous Lorentz transformation.

The set of all holonomy operators at a vertex $x$ is obviously a subgroup of $O(1,n-1)$ (or $SO(1,n-1)^+$ in case of a spin connection), and since we assume $G$ to be connected there is but one holonomy group up to canonical isomorphism. We call it \emph{the holonomy group of $G$}. We say that $G$ is flat when its holonomy group is trivial.

The discrete connection will be said to be \emph{a spin connection} if it is  Clifford, metric, spin and orientation preserving.

\subsection{The split Dirac structure as a spectral spacetime}

We might be expecting that the split Dirac structure is a spectral spacetime exactly when the discrete connection is spin. It turns out to be \emph{almost} true.

\begin{theorem}\label{th6} The split  Dirac structure over the  graph  $G$ is a (possibly non-reconstructible) antilorentzian spectral spacetime if and only if
\begin{enumerate}
\item The elements $\gamma_e^\pm$ are Krein self-adjoint (hence odd),  satisfy 

\be
J_{e^+}\gamma_e^+(J_{e^+})^{-1}=-h_e^+\gamma_e^-h_e^-\label{eq83}
\ee

for all $e$, and do not vanish identically on the left and right component of $S_{e^\pm}$.
\item The discrete connection is metric, spin and orientation preserving.
%\item The holonomy group fixes an odd, $c$-real, and Krein-positive element.
%\item $\Gamma_{e^+}=h_e^+\Gamma_{e^-}h_e^-$.
\end{enumerate}
When this is so, the KO signs are $\tilde\epsilon=\epsilon,\tilde\epsilon''=\epsilon''$ and  positive orientation 1-forms are of the form:

$$\beta F(e,\pm)=  \Gamma_e^\pm  h_e^\pm F(e,\mp)$$
where $\Gamma_e^\pm$ is Krein-positive and satisfy 

\be
h_e^\mp\Gamma_e^\pm h_e^\pm=-J_{e^\mp}\Gamma_e^\mp J_{e^\mp}^{-1}\label{condforms}
\ee
\end{theorem}

Before embarking into the proof, let us observe that the most obvious solution to (\ref{eq83})  is to take $\gamma_e^\pm$ to be equal to $\rho(v_e^\pm)$ where $v_e^\pm$ is a  non-zero   vector in $X_{e^\pm}$.  In this case we say that the split Dirac structure is \emph{vectorial}. If in addition  the family $\{\gamma_e^\pm|e^\pm=v\}$  generates the space $X_v$ at each vertex $v$, we say that it is \emph{complete}. If the split Dirac structure is vectorial and complete  we see that (\ref{eq83}) forces   the discrete connection to be a Clifford connection, and will thus be a spin connection. Condition (\ref{eq83}) can then be seen as a remnant of the missing Clifford property of the discrete connection.  %This will become clearer in section \ref{secdiscr} when we see the relation with the discretization of the usual Dirac operator.

Moreover, condition (\ref{condforms}) can always be met by choosing $\Gamma_e^-$ to be a $J$-imaginary (i.e. $c$-real) Krein-positive operator (for instance a future-directed timelike vector), and let $\Gamma_e^+$ be defined by the LHS of (\ref{condforms}). It will then be Krein-positive since the $h_e^\pm$ are Krein-unitary.

\begin{demo}
We begin by observing that (\ref{defkrein}) defines a sesquilinear form iff $h_e^\pm$ is Krein unitary, \emph{hence the discrete connection must be metric}.

Now we ask that $\chi^\times=-\chi$.

\bea
(\chi F,G)&=&\sum_{e\in E}\left((\chi_{e^+} F(e,+),h_e^+G(e,-))+( \chi_{e^-}F(e,-),h_e^-G(e,+))\right)\cr
&=&-\sum_{e\in E}\left((F(e,+),\chi_{e^+}h_e^+G(e,-))+(F(e,-), \chi_{e^-}h_e^-G(e,+))\right)\nonumber
\eea

\emph{Hence we see that  $\chi^\times=-\chi$ iff
 $\chi_{e^\pm} h_e^\pm=h_e^\pm \chi_{e^\mp}$, that is, if the discrete connection is orientation preserving.} We suppose this from now on.

Next we check under what condition  $D$ anticommutes with $\chi$:

\bea
D(\chi F)(e,\pm)&=&{\mp  1\over \delta_e}\gamma_e^\pm h_e^\pm \chi_{e^\mp} F(e,\mp)\cr
&=&{\mp  1\over \delta_e}\gamma_e^\pm \chi_{e^\pm}h_e^\pm  F(e,\mp)
\eea
and
\be
\chi (D F)(e,\pm)={\mp  1\over \delta_e}\chi_{e^\pm}\gamma_e^\pm h_e^\pm F(e,\mp)\nonumber
\ee

Thus we see that \emph{$D$ anticommutes with $\chi$ iff $\gamma_e^\pm$ is odd.}

We check whether $D=D^\times$. We have:
 
\bea
(DF,G)&=&\sum_{e\in E}\left((DF(e,+), h_e^+G(e,-)+(DF(e,-), h_e^-G(e,+))\right)\cr
&=&\sum_{e\in E}\left(({-1\over \delta_e}\gamma_e^+h_e^+F(e,-),h_e^+G(e,-))+({1\over \delta_e} \gamma_e^-h_e^-F(e,+),h_e^-G(e,+))\right)\cr
&=&\sum_{e\in E}\left(({-1\over \delta_e} F(e,-),h_e^- (\gamma_e^+)^\times h_e^+G(e,-))+({1\over \delta_e}  F(e,+),h_e^+ (\gamma_e^-)^\times h_e^-G(e,+))\right)\nonumber
\eea

and
\bea
(F,DG)&=&\sum_{e\in E}\left((F(e,+), h_e^+DG(e,-)+(F(e,-), h_e^-DG(e,+))\right)\cr
&=&\sum_{e\in E}\left({1\over \delta_e}(F(e,+), h_e^+\gamma_e^-h_e^-G(e,+)-{1\over\delta_e}(F(e,-), h_e^-\gamma_e^+h_e^+G(e,-))\right)\nonumber
\eea

\emph{Hence we see that $D=D^\times$ iff $\gamma_e^\pm$ is Krein self-adjoint}, which we suppose from now on.

Now let us explore the properties of $J$. We want $J^2=\tilde\epsilon=\pm 1$.

\bea
(J^2F)(e,\pm)&=& h_e^\pm J_{e^\mp} (JF)(e,\mp)\cr
\Leftrightarrow \tilde\epsilon F(e,\pm) &=&h_e^\pm J_{e^\mp} h_e^\mp J_{e^\pm} F(e,\pm)\nonumber 
%\Leftrightarrow \tilde\epsilon(J_{e^\pm})^{-1}&=& h_e^\pm J_{e^\mp} h_e^\mp
\eea
Thus we obtain $\tilde\epsilon(J_{e^\pm})^{-1}= h_e^\pm J_{e^\mp} h_e^\mp$. \emph{Hence we see that $J^2=\tilde\epsilon$ iff }
% \bea
% \Leftrightarrow \epsilon\tilde\epsilon J_{e^\pm}&=&h_e^\pm J_{e^\mp} h_e^\mp\cr
% \eea

\be
h_e^\pm J_{e^\mp} h_e^\mp=\epsilon\tilde\epsilon J_{e^\pm}\label{ref54}
\ee

Now we compute $DJ$ and $JD$. We find that

% \bea
% (DJ)F(e,\pm)&=&{\mp 1\over\delta_e}\gamma_e^\pm h_e^\pm JF(e,\mp)\cr
% &=&{ \mp 1\over\delta_e}\gamma_e^\pm h_e^\pm h_e^\mp J_{e^\pm} F(e,\pm)\cr
%  &=&{ \mp 1\over\delta_e}\gamma_e^\pm  J_{e^\pm} F(e,\pm) 
% \eea

$$(DJ)F(e,\pm)={ \mp 1\over\delta_e}\gamma_e^\pm  J_{e^\pm} F(e,\pm) $$

and

$$(JD)F(e,\pm)={\pm 1\over\delta_e} \epsilon\tilde\epsilon J_{e^\pm} h_e^\pm \gamma_e^\mp h_e^\mp F(e,\pm)$$

% \bea
% (JD)F(e,\pm)&=& h_e^\pm J_{e^\mp} DF(e,\mp)\cr
% &=&{\pm 1\over\delta_e} h_e^\pm J_{e^\mp}\gamma_e^\mp h_e^\mp F(e,\pm)\cr
% &=&{\pm 1\over\delta_e} \epsilon\tilde\epsilon J_{e^\pm} h_e^\pm \gamma_e^\mp h_e^\mp F(e,\pm)
% \eea

\emph{We see that $JD=DJ$ iff}
 
\be
J_{e^\pm} \gamma_e^\pm (J_{e^\pm})^{-1}=-\epsilon\tilde\epsilon h_e^\pm \gamma_e^\mp h_e^\mp\label{eq57}
\ee
 which we suppose from now on. We now compare $J\chi$ and $\chi J$:

\bea
(J\chi) F(e,\pm)&=&  h_e^\pm J_{e^\mp} \chi_{e^\mp} F(e,\mp)\cr
&=&\epsilon'' h_e^\pm \chi_{e^\mp}J_{e^\mp} F(e,\mp)\nonumber
\eea

and 

\bea
(\chi J)F(e,\pm)&=& \chi_{e^\pm}JF(e,\pm)\cr
&=&\chi_{e^\pm}h_e^\pm J_{e^\mp} F(e,\mp)\cr
&=&h_e^\pm \chi_{e^\mp}J_{e^\mp} F(e,\mp)\nonumber
\eea

where in the last line we have used the fact that the discrete connection preserves spacetime orientation. Hence we obtain $\tilde \epsilon''=\epsilon''$.

Now we compute $J^\times$.

\bea
(JF,G)&=&\sum_{e\in E}\big((JF(e,+),h_e^+G(e,-))+(JF(e,-),h_e^-G(e,+)\big)\cr
&=&\sum_{e\in E}\big((h_e^+J_{e^-}F(e,-),h_e^+G(e,-))+(h_e^-J_{e^+}F(e,+),h_e^-G(e,+)\big)\cr
&=&\sum_{e\in E}\big((J_{e^-}F(e,-), G(e,-))+(J_{e^+}F(e,+),G(e,+)\big) \nonumber
\eea

and

\bea
(F,JG)&=&\sum_{e\in E}\big((F(e,+),h_e^+JG(e,-))+(F(e,-),h_e^-JG(e,+)\big)\cr
&=&\sum_{e\in E}\big((F(e,+),h_e^+h_e^-J_{e^+}G(e,+))+(F(e,-),h_e^-h_e^+J_{e^-}G(e,-)\big)\cr
&=&\sum_{e\in E}\big((F(e,+),J_{e^+}G(e,+))+(F(e,-),J_{e^-}G(e,-)\big)\nonumber
\eea

so that

\bea
\overline{(F,JG)}&=&\sum_{e\in E}\big(((J_{e^+})^\times F(e,+),G(e,+))+((J_{e^-})^\times F(e,-),G(e,-)\big)\cr
&=&\epsilon\kappa\sum_{e\in E}\big((J_{e^+}F(e,+),G(e,+))+(J_{e^-}F(e,-),G(e,-)\big)\nonumber
\eea

Now we want $J^\times=\tilde\epsilon\tilde\kappa J$. Using the fact that we required $\tilde\kappa=\kappa=-1$, we obtain  $\epsilon =\tilde\epsilon$. Therefore we get from (\ref{ref54}) that the discrete connection is spin preserving.

Let us now consider the orientation forms.  First we find that

\bea
[D,\pi(a)]F(e,\pm)&=& {\pm 1\over\delta_e}\left(a(e^\pm)\gamma_e^\pm h_e^\pm-\gamma_e^\pm h_e^\pm a(e^\mp)\right)F(e,\mp)\cr
&:=& \pm{\delta a_e^\pm\over \delta_e} F(e,\mp)\label{eq61}
\eea

% \bea
% [D,\pi(a)]F(e,\pm)&=&(D(\pi(a)F)(e,\pm)-(\pi(a)(DF))(e,\pm)\cr
% &=&{\mp 1\over  \delta_e}\gamma_e^\pm h_e^\pm(\pi(a)F)(e,\mp)-a(e^\pm)(DF)(e,\pm)\cr
% &=&{\mp 1\over  \delta_e}\gamma_e^\pm h_e^\pm a(e^\mp)F(e,\mp)+{\pm 1\over \delta_e}a(e^\pm)\gamma_e^\pm h_e^\pm F(e,\mp)\cr
% &=& \pm{\delta a_e^\pm\over \delta_e} F(e,\mp)\label{eq61}
% \eea

%where we have introduced $\delta a_e^\pm=a(e^\pm)\gamma_e^\pm h_e^\pm-\gamma_e^\pm h_e^\pm a(e^\mp)$. 

Now let $f=(v,v')$ be a particular edge and introduce the scalar functions $\delta_v$ and $\delta_{v'}$ which take the value zero everywhere except on the named vertices. Then we have

\bea
\pi(\delta_{v'})([D,\pi(\delta_{v})]F)(e,\pm)&=&\delta_{v'}(e^\pm)[D,\pi(\delta_{v})]F(e,\pm)\cr
&=&{\pm 1\over\delta_e}\delta_{v'}(e^\pm)((\delta_{v}(e^\pm)-\delta_{v}(e^\mp))\gamma_e^\pm h_e^\pm F(e,\mp)\nonumber
\eea

One sees that this expression vanishes unless $\{e^+;e^-\}=\{v;v'\}$, and since we have only one edge joining two vertices, this means that $e=f$, hence $e^-=v$ and $e^+=v'$. To sum up, we have that  $\pi(\delta_{v'})[D,\pi(\delta_{v})]F$ vanishes everywhere on the split graph $\tilde E$, except at $(f,+)$ where its value is $-{1\over\delta_e}\gamma_e^+h_e^+F(f,-)$. We can of course make a similar reasoning with $\pi(\delta_{v})([D,\pi(\delta_{v'})]$ and obtain a 1-form which vanishes everywhere except at $(f,-)$ where its value is ${1\over\delta_e}\gamma_e^-h_e^-F(f,+)$. Summing these two 1-forms we obtain an element of $\Omega^1_D({\cal A}^0,\pi)$ which has the following structure:

\be
\matrix{ & S_{e^-}^L& S_{e^-}^R & S_{e^+}^L & S_{e^+}^R \cr 
S_{e^-}^L& 0 & 0 &0 &A^-\cr
S_{e^-}^R& 0 & 0 &B^- &0\cr
S_{e^+}^L& 0 & A^+ &0 &0\cr
S_{e^+}^R& B^+ & 0 &0 &0
}
\label{strucform}
\ee

where we have decomposed $S_{e^\pm}$ into its left and right components, and $A^\pm,B^\pm$ stand for  possibly non-zero elements. The    block $\pmatrix{0&A^\pm\cr B^\pm&0}$ is $\gamma_e^\pm h_e^\pm$. The algebra ${\cal A}^0$ acts with block-diagonal matrices like

$$\matrix{ & S_{e^-}^L& S_{e^-}^R & S_{e^+}^L & S_{e^+}^R \cr 
S_{e^-}^L& a_L^- & 0 &0 & 0\cr
S_{e^-}^R& 0 & a_R^- & 0 &0\cr
S_{e^+}^L& 0 & 0 & a_L^+ &0\cr
S_{e^+}^R& 0 & 0 &0 &a_R^+
}
$$

and since $\Omega^1_D$ is an ${\cal A}^0$-bimodule, we can act on both sides. Thus we see that  a general 1-form is just any operator with the same structure as (\ref{strucform}). In particular, if a block among $A^\pm,B^\pm$ vanishes, we see that no 1-form is invertible, hence  no 1-form can be Krein-positive. Hence we have a restriction on $\gamma_e^\pm$: we see that \emph{$\gamma_e^\pm$ must not vanish identically on $S_{e^\pm}^L$ or $S_{e^\pm}^R$}. We make this hypothesis from now on.

Thus a general 1-form $\beta$ acts as

\be
\beta F(e,\pm)=\beta_e^\pm F(e,\mp)\label{eq63}
\ee

where $\beta_e^\pm$ is any  linear operator from  $S_{e^\mp}$ to $S_{e^\pm}$ which exchanges chiralities, and can thus be written

\be
\beta_e^\pm=\Gamma_e^\pm h_e^\pm,\mbox{ with }\chi_{e^\pm}\Gamma_e^\pm=-\Gamma_e^\pm\chi_{e^\pm}
\ee

We now look for conditions on the $\Gamma_e^\pm$ for $\beta$ to be Krein self-adjoint and imaginary. For this we first remark that

\be
(\beta^\times F)(e,\pm)=h_e^\pm (\beta_e^\pm)^\times h_e^\pm F(e,\mp)
\ee

% Indeed, we have:

% \bea
% (\beta F,G)&=&\sum_e \left((\beta F(e,+),h_e^+ G(e,-))+(\beta F(e,-),h_e^- G(e,+))\right)\cr
% &=&\sum_e \left((\beta_e^+F(e,-),h_e^+G(e,-))+(\beta_e^-F(e,+),h_e^-G(e,+))\right)\cr
% &=&\sum_e \left((F(e,-),(\beta_e^+)^\times h_e^+G(e,-))+(F(e,+),(\beta_e^-)^\times h_e^-G(e,+))\right)\cr
% &=&\sum_e \left((F(e,-),h_e^-h_e^+(\beta_e^+)^\times h_e^+G(e,-))+(F(e,+),h_e^+h_e^-(\beta_e^-)^\times h_e^-G(e,+))\right)\cr
% &=&(F,\beta^\times G)
% \eea

So that 

\be
\beta^\times=\beta\Leftrightarrow (\beta_e^\pm)^\times=h_e^\mp\beta_e^\pm h_e^\mp\Leftrightarrow (\Gamma_e^\pm)^\times=\Gamma_e^\pm\label{cond65}
\ee

Moreover, we also have:

\be
J\beta J^{-1}=-\beta\Longleftrightarrow \beta_e^\pm=- h_e^\pm J_{e^\mp}\beta_e^\mp J_{e^\pm}^{-1} h_e^\pm\Leftrightarrow h_e^\mp\Gamma_e^\pm h_e^\pm=-J_{e^\mp}\Gamma_e^\mp J_{e^\mp}^{-1}\label{cond66}
\ee

% {\small Check:
% \bea
% J\beta J^{-1}F(e,\pm)&=&\tilde \epsilon J\beta JF(e,\pm)\cr
% &=&\tilde \epsilon h_e^\pm J_{e^\mp} (\beta J F)(e,\mp)\cr
% &=&\tilde \epsilon h_e^\pm J_{e^\mp} \beta_e^\mp (JF)(e,\pm)\cr
% &=&\tilde \epsilon h_e^\pm J_{e^\mp} \beta_e^\mp h_e^\pm J_{e^\mp} F(e,\mp)\cr
% &=&\epsilon h_e^\pm J_{e^\mp} \beta_e^\mp J_{e^\pm} h_e^\pm F(e,\mp)
% \eea
% }
 
Finally it is easy to see that \emph{$\bra .,.\ket_\beta$ is positive definite iff $\Gamma_e^\pm$ is Krein positive.}
\end{demo}

\subsection{The reconstructibility of the split Dirac structure}
Now let us investigate the reconstructibility condition $\pi({\cal A})^{*_\beta}=\pi({\cal A})$. 

 We first note that

\be
\beta^{-1}F(e,\pm)=(\beta_e^\mp)^{-1}F(e,\mp)
\ee

and 

\be
\pi(a)^\times F(e,\pm)=h_e^\pm a(e^\mp)^\times h_e^\mp F(e,\pm)\label{eq73}
\ee

hence

\bea
\beta\pi(a)^\times\beta^{-1} F(e,\pm)&=&\beta_e^\pm h_e^{\mp}a(e^\pm)^\times h_e^\pm(\beta_e^\pm)^{-1}F(e,\pm)\cr
&=&\Gamma_e^\pm a(e^\pm)^\times (\Gamma_e^\pm)^{-1}F(e,\pm)\label{eq74}
\eea

which we want to be of the form $b(e^\pm)F(e,\pm)$ for some $b\in{\cal A}$. This creates constraints on the operators $\Gamma_e^\pm$. Indeed, if two edges share a vertex $v$, for instance if $e^+=f^+$, then $Ad_{\Gamma_e^+}=Ad_{\Gamma_f^+}$, hence $\Gamma_e^+(\Gamma_f^+)^{-1}$ commutes with $\CC l(X_v)^0$ and thus belong to $\CC\oplus \CC\chi_v$. If $n=2$ this constraint is always satisfied since $\Gamma_e^+(\Gamma_f^+)^{-1}\in \CC l(V)^0=\CC\oplus\CC\chi_v$ in this case. Hence we can already conclude that the split Dirac structure is always reconstructible when $n=2$.

We now suppose that $n>2$. We will able to conclude only when the connection is Clifford, so we suppose this property holds in what follows. For each vertex $v$ we can choose arbirtrarily an edge $f$ such that $f^+$ or $f^-$ is equal to $v$, and set  $\Gamma_v:=\Gamma_f^\pm$, according to the case. Then for every edge $e$ we will have the relation

\be
\Gamma_e^\pm=\lambda_e^\pm \Gamma_{e^\pm}
\ee

where $\Gamma_{e^\pm}$ is Krein positive and $\lambda_e^\pm=x_e^\pm+y_e^\pm \chi_{e^\pm}$, with $x_e^\pm,y_e^\pm\in\RR$. %Using the fact that $\Gamma_e^\pm$ is Krein positive we obtain that $((x_e^\pm+y_e^\pm \chi_{e^\pm})\Gamma_{e^\pm}\psi,\psi)>0$ for every non-zero $\psi\in S_{e^\pm}$. From this we see that $x_e^\pm-y_e^\pm\chi_{e^\pm}$ is a positive operator for the scalar product $\bra .,.\ket_{\Gamma_{e^\pm}}$, hence $x_e^\pm>|y_e^\pm|\ge 0$.

% , which we can take to be real, and $\Gamma_{e^\pm}^\times=\Gamma_{e^\pm}$, thanks to conditions (\ref{cond65}). Moreover we must have $\lambda_e^\pm\Gamma_{e^\pm}$ Krein positive. However, thanks to condition (\ref{cond66}) we have

% \be
% \lambda_e^\pm h_e^\mp \Gamma_{e^\pm} h_e^\pm=-\lambda_e^\mp J_{e^\mp} \Gamma_{e^\mp} (J_{e^\mp})^{-1}\label{cond70}
% \ee

% This shows that if the vertex $v$ is shared by two edges, for instance if $v=e^+=f^-$ then $\lambda_e^+$ and $\lambda_f^-$ must have the same sign. Hence we can arrange, by putting this sign into $\Gamma_v$, that $\lambda_e^\pm$ is always positive. 
% {\small details :
% \bea
% \sum_e\big((F(e,+),h_e^+\beta^{-1}F(e,-))+(F(e,-),h_e^-\beta^{-1}F(e,+))\big)&\ge &0,\forall F\cr
% \sum_e\big((F(e,+),h_e^+(\beta_e^+)^{-1}F(e,+))+(F(e,-),h_e^-(\beta_e^-)^{-1}F(e,-))\big)&\ge &0\cr
% \sum_e\big((F(e,+),(\Gamma_e^+)^{-1}F(e,+))+(F(e,-),(\Gamma_e^-)^{-1}F(e,-))\big)&\ge &0\cr
% \eea

% }
%  

% \bea
% r_e(F(e,\pm),h_e^\pm\beta_e^\mp F(e,\mp))&=& \lambda_e^\pm  r_e((F(e,\pm),\Gamma_{e^\pm} F(e,\pm))\cr
% &>&0\mbox{ for all }e\in E
% \eea

% Let us sum up the constraints which have appeared on $\beta$:

%    

% \be
% \beta_e^\pm=\lambda_e^\pm \Gamma_{e^\pm}h_e^\pm
% \ee

% with $\lambda_e^\pm, \Gamma_{e^\pm}$ such that $\lambda_e^\pm>0$,  and  $\Gamma_{e^\pm}$ is Krein-positive (hence odd and  Krein-self-adjoint), and satisfies   condition (\ref{cond70}).

Now let $H_e^\pm : \CC l(X_{e^\mp})\rightarrow \CC l(X_{e^\pm})$ be the map $a\mapsto h_e^\pm a h_e^\mp$   and $\tilde H_e^\pm=c_{e^\pm}\circ H_e^\pm$. We can imagine $H_e^\pm$ to be the canonical extension to the Clifford algebras of the parallel transport of the Levi-Civita connection, and we note that  thanks to $J_{e^\pm} h_e^\pm=h_e^\pm J_{e^\mp}$, we have $c_{e^\pm}\circ H_e^\pm=H_e^\pm\circ c_{e^\mp}$. Then we can give constraint (\ref{condforms}) the nicer form:

\be
\Gamma_{e^\pm}=\mu_e^\pm \tilde H_e^\pm(\Gamma_{e^\mp})\label{eq77}
\ee

where $\mu_e^\pm=(\lambda_e^\pm)^{-1}\tilde H_e^\pm(\lambda_e^\mp)$.  

Now let ${\cal P}=(v_1,\ldots,v_k)$ be a path from $v=v_1$ to $v'=v_k$. If $(v_i,v_{i+1})$ is an edge $e_i$ of $G$ we set $\sigma_i=+$, if $(v_{i+1},v_i)$ is an edge we set $\sigma_i=-$. Using (\ref{eq77}) repeatedly along ${\cal P}$ we obtain

\be
\Gamma_{v'}=\mu_{\cal P} \tilde H_{\cal P} (\Gamma_{v})\label{holconstr}
\ee

where 

\be
\mu_{\cal P}=\prod_{i=1}^{k-1}\tilde H_{{\cal P}_{i+1}}(\mu_{e_i}^{\sigma_i}),\mbox{ with } \tilde H_{{\cal P}_i}=\prod_{j=i}^{k-1}\tilde H_{e_j}^{\sigma_j},\mbox{ and }\tilde H_{\cal P}=\tilde H_{{\cal P}_1}
\ee

%From the remark above we can see that $\mu_{\cal P}$ has the same sign as the relative signs of the $\lambda$'s at the extremities of ${\cal P}$.

(Hence $\tilde H_{{\cal P}_i}$ is the parallel transport along the subpath $(v_i,\ldots,v_k)$, composed with the real structure $c_{v_k}$ if the length of the subpath is odd.)

We note that if ${\cal P}_1$ is a path from $v_1$ to $v$ and ${\cal P}_2$ is a path from $v$ to $v_2$ we have

\be
\mu_{{\cal P}_1*{\cal P}_2}=\mu_{{\cal P}_2}\tilde H_{{\cal P}_2}(\mu_{{\cal P}_1})\label{prodmu}
\ee 

Now suppose ${\cal P}$ is a closed path from $v$ to $v$. We write $\tilde H_{\cal P}={\rm Hol}_{\cal P}\circ c_v$, where ${\rm Hol}_{\cal P}=Ad_{{\rm hol}_{\cal P}}$ is the holonomy around ${\cal P}$. Here ${{\rm hol}_{\cal P}}$, which is the product of the $h_e^\pm$ along the path, is Krein-unitary and commutes with $\chi_v$. Hence ${{\rm hol}_{\cal P}}$ commutes with $\mu_{\cal P}=x+y\chi_v$. For each $\psi\in S_v$ we then have

\bea
(\psi,\Gamma_v \psi)&=&(\psi,\mu_{\cal P}{{\rm hol}_{\cal P}}c(\Gamma_v){{\rm hol}_{\cal P}}^{-1}\psi)\cr
&=&({{\rm hol}_{\cal P}}^{-1}\psi,\mu_{\cal P}c(\Gamma_v){{\rm hol}_{\cal P}}^{-1}\psi)\cr
&=&-({{\rm hol}_{\cal P}}^{-1}\psi,(x+y\chi_v)J_v\Gamma_v J_v^{-1}{{\rm hol}_{\cal P}}^{-1}\psi)\cr
&=&-({{\rm hol}_{\cal P}}^{-1}\psi,J_v(x-\epsilon'' y\chi_v)\Gamma_v J_v^{-1}{{\rm hol}_{\cal P}}^{-1}\psi)\cr
&=&(J_v^{-1}{{\rm hol}_{\cal P}}^{-1}\psi,(x-\epsilon'' y\chi_v)\Gamma_v J_v^{-1}{{\rm hol}_{\cal P}}^{-1}\psi)\cr
&=&((x+\epsilon'' y\chi_v)J_v^{-1}{{\rm hol}_{\cal P}}^{-1}\psi,\Gamma_v J_v^{-1}{{\rm hol}_{\cal P}}^{-1}\psi)\cr
&=&\bra (x+\epsilon'' y\chi_v)\psi',\psi'\ket_{\Gamma_v},\mbox{ with }\psi'=J_v^{-1}{{\rm hol}_{\cal P}}^{-1}\psi
\eea

Hence we see that $x+\epsilon''y\chi_v$ is a positive operator, which entails that $x>|y|$.

Turning twice around ${\cal P}$ we obtain from (\ref{prodmu}) 

\be
\Gamma_v=\mu_{\cal P}c(\mu_{\cal P}) {\rm Hol}_{\cal P}^{2}(\Gamma_v)
\ee

if ${\cal P}$ has an odd length, and

\be
\Gamma_v=\mu_{\cal P}^{2} {\rm Hol}_{\cal P}^{2}(\Gamma_v)
\ee

if ${\cal P}$ is of even length. In the second case taking the determinant  yields $x^2-y^2=\pm 1$, hence $x^2-y^2=1$ since $x>|y|$. In the first case we obtain that $\det(\mu_{\cal P}c(\mu_{\cal P}))=1$. If $\epsilon''=-1$, $c(\mu_{\cal P})=\mu_{\cal P}$ and we obtain $x^2-y^2=1$ again. If $\epsilon''=1$, $\mu_{\cal P}c(\mu_{\cal P})=x^2-y^2$ so the conclusion still holds. We will write $\mu_{\cal P}=\cosh t+\chi_v\sinh t$ in the sequel.

{\bf  Let us first suppose that $\epsilon''=-1$.}

We consider a pseudo-orthonormal basis $(e_i)_{1\le i\le n}$ in $X_v$ and the associated pseudo-orthonormal basis $(e_I)_{I\subset\{1;\ldots;n\}}$ in $\CC l(X_v)$. We decompose $\CC l(X_v)^1$ as an orthogonal sum $\CC l(X_v)^1=X_v^1\oplus \ldots\oplus X_v^{n-1}$ as in \cite{part1}, proposition \ref{propo4}, and we note that $X_v^1=X_v$ and $X_v^{n-1}=\chi_v X_v$. Hence, writing $L$ for the sum of the middle terms, we have the $(.,.)_c$-orthogonal decomposition $\CC l(X_v)=X_v\oplus L\oplus \chi_v X_v$. Let us write $\Gamma_v=u+l+\chi_v w$, according to this decomposition.

From (\ref{holconstr}) we obtain, using  the conservation of the orientation:

\be
(\cosh t-\chi_v\sinh t)(u+\ldots+\chi_v w)=(\Lambda u+\ldots+\chi_v\Lambda w)\label{eq80}
\ee

where $\Lambda$ is some proper orthochronous Lorentz transformation and the dots stand  for some elements in $L$. Note that (\ref{eq80}) holds for paths of both odd and even lengths thanks to the properties $c(u)=u$ and $c(\chi_v w)=\chi_v w$. Using proposition \ref{propo4} in \cite{part1} again we obtain

\be
\left\{
\matrix{\Lambda u&=& u\cosh t-w\sinh t\cr \Lambda w&=&-u\sinh t+w\cosh t}
\right.\label{eq82}
\ee

Now the idea is that $\Lambda$ must be a boost to have this action on $u$ and $w$, which amounts to say that $g(u,u)+g(w,w)=0$ and $g(u,w)=0$. To see this, we compute $(\Lambda u)^2+(\Lambda w)^2$ and $g(u,w)$ from (\ref{eq82}). We obtain

\be
\left\{
\matrix{u^2+w^2&=&(\Lambda u)^2+(\Lambda w)^2&=& (u^2+w^2)\cosh(2t)-2\sinh(2t)g(u,w) \cr 
g(u,w)&=&g(\Lambda u,\Lambda w)&=&-{1\over 2}\sinh(2t)(u^2+w^2)+\cosh(2t)g(u,w)\cr}
\right.  \nonumber
\ee
Thus
\be
\left\{
\matrix{  (\cosh(2t)-1)(u^2+w^2)-2\sinh(2t)g(u,w)&=&0 \cr 
 -{1\over 2}\sinh(2t)(u^2+w^2)+(\cosh(2t)-1)g(u,w)&=&0\cr}
\right.  \nonumber
\ee

The determinant of the system vanishes only if $t=0$. Suppose that $t\not=0$. Then $u^2+w^2=g(u,w)=0$, hence $(u+w)^2=0$. Thus $u+w$ is lightlike. Moreover we know from \cite{part1}, lemma  \ref{lem9},  that $u+\chi_v w+l$  Krein-positive implies that $u$ is future-directed timelike. For every $\epsilon\in]0;2[$, let $x_\epsilon=u-(1-\epsilon)\chi_vw$. Since  $u\pm(1-\epsilon)w$ is future-directed timelike, we conclude from   \cite{part1}, lemma \ref{lem8}   that $x_\epsilon$  is Krein-positive. Since $\Gamma_v$ is also Krein-positive, this entails that $ x_\epsilon\Gamma_v$ is a positive operator for the scalar product $\bra .,.\ket_{\Gamma_v}$, and in particular that $\tau_n(x_\epsilon\Gamma_v)>0$, where $\tau_n$ is the normalized trace on $\CC l(X_v)$.  However, we have:

\bea
\tau_n\left(x_\epsilon(u+\chi_v w+l)\right)&=&\tau_n\left(\epsilon(-g(w,w)+\chi_v wu)+x_\epsilon l\right)\cr
&=&-\epsilon g(w,w) 
\eea

Hence we see taking the limit $\epsilon\rightarrow 0$ that the positive semi-definite operator $(u-\chi_v w)(u+\chi_v w+l)$ has zero trace and is therefore equal to zero. This contradicts the fact that $u+\chi_v w+l$ is  Krein-positive, hence invertible.

We can then conclude that $t=0$, hence $\mu_{\cal P}=1$. Thus we have $\Gamma_v=\tilde H_{\cal P}(\Gamma_v)$ for each loop ${\cal P}$ based at $v$.

{\bf Suppose now that $\epsilon''=1$.} In this case the exact same conclusion holds if ${\cal P}$ is a closed path at $v$ of even length. Moreover, if ${\cal P}_1$ and ${\cal P}_2$  are paths of odd length, then $\mu_{{\cal P}_1*{\cal P}_2}=\mu_{{\cal P}_1}c(\mu_{{\cal P}_2})=\mu_{{\cal P}_1}\mu_{{\cal P}_2}^{-1}=1$ since ${\cal P}_1*{\cal P}_2$ is even. Hence $\mu_{\cal P}$ only depends on the parity of ${\cal P}$: it is equal to $1$ if ${\cal P}$ is even, and to $\cosh t_v+\chi_v\sinh t_v$ if ${\cal P}$ is odd, where $t_v$ only depends on $v$. Now let $\tilde \Gamma_v:=x_v\Gamma_v$, where $x_v=\cosh(t_v/2)-\sinh(t_v/2)\chi_v$. If ${\cal P}$ is even we have $\tilde H_{\cal P}(\tilde \Gamma_v)={\rm Hol}_{\cal P}(x_v\Gamma_v)=x_v\Gamma_v=\tilde\Gamma_v$, and if ${\cal P}$ is odd we have $\tilde H_{\cal P}(\tilde \Gamma_v)=c_v(x_v)\tilde H_{\cal P}(\Gamma_v)=x_v^{-1}\mu_{\cal P}^{-1}\Gamma_v$. Hovewer $x_v^{-1}\mu_{\cal P}^{-1}=(\cosh(t_v/2)+\sinh(t_v/2)\chi_v)(\cosh( t_v)-\sinh(t_v)\chi_v)=\cosh(t_v/2)-\sinh(t_v/2)\chi_v=x_v$. Hence $\tilde H_{\cal P}(\tilde \Gamma_v)=\tilde\Gamma_v$ in every case.

Thus, we have obtained for $\epsilon''=-1$ and $\epsilon''=1$ alike the existence of a Krein-positive element $\tilde\Gamma_v$ which satisfies $\tilde H_{\cal P}(\tilde\Gamma_v)=\tilde\Gamma_v$ for every closed path based at $v$. Let us write $\tilde \Gamma_v=u_v+\ldots$ using the orthogonal decomposition $\CC l(V)^1=V^1\oplus\ldots\oplus V^{n-1}$, which is preserved by the holonomy. Since $c_v(u_v)=u_v$ we obtain that 

\be
{\rm Hol}_{\cal P}(u_v)=u_v
\ee

for every closed path at $v$. Hence for every vertex $w$ we can define 

\be
u_{w}:=H_{\cal Q}(u_v)\label{cc}
\ee

where ${\cal Q}$ is any path from $v$ to $w$. Of course  this will not depend on the chosen path. Moreover we know from \cite{part1}, lemma \ref{lem8}, that $u_{v}$, hence $u_{w}$ will be future-directed timelike. On a manifold, a vector field which is parallel transported to itself as in (\ref{cc}) is said to be parallel, or covariantly constant, since its covariant derivative must vanish in every direction. This justifies the following definition.

\begin{definition}
A vector field on $G$, that is a section $v\mapsto u_v$ of the discrete vector bundle $\bigcup_{v\in G}X_v$, is said to be \emph{parallel}, or \emph{covariantly constant},  if for all $w\in G$ and for all path ${\cal Q}$ from $v$ to $w$, one has

$$u_{w}=H_{\cal Q}(u_v)$$
\end{definition}

 We have therefore obtained a covariantly constant future-directed timelike vector field on $G$ under the assumption that the split Dirac structure were reconstructible.

Conversely, if $v\mapsto u_v$ is a covariantly constant future-directed timelike vector field on $G$, then let $\beta$ be defined by $\beta_e^\pm=\Gamma_e^\pm h_e^\pm$ with $\Gamma_e^\pm:=u_{e^\pm}$.

With this definition $\beta$ is a positive orientation 1-form since it satisfies $h_e^\pm\Gamma_e^\mp h_e^\mp=H_e^\pm (u_{e^\mp})=u_{e^\pm}=c(u_{e^\pm})=-J_{e^\pm}u_{e^\pm}J_{e^\pm}^{-1}$. Moreover from equation (\ref{eq74}) it satisfies

\bea
\beta(\pi(a)^\times)\beta^{-1}F(e,\pm)&=&u_{e^\pm}a(e^\pm)^\times u_{e^\pm}^{-1}F(e,\pm)\cr
&=&b(e^\pm)F(e,\pm)
\eea

where $b(e^\pm)\in \CC l(X_{e^\pm})^0$. Thus the split Dirac structure is reconstructible.

All in all we have proven the following theorem.

\begin{theorem}\label{th7} Suppose the hypotheses of theorem \ref{th6} are met. Then if $n=2$,  the split  Dirac structure over  $G$ is reconstructible. If $n>2$ and the discrete connection is Clifford, then the split Dirac structure  is reconstructible iff there exists a  future-directed timelike and parallel vector field on $G$.
%
% if and only if the holonomy group fixes an odd, $c$-real, and Krein-positive element. In this case the orientation forms which satisfy the reconstructibility condition are 
%$$\beta F(e,\pm)=k_e \Gamma_{e^\pm} h_e^\pm F(e,\mp)$$
%with $k_e>0$ for all $e\in E$, and $\Gamma_{e^\pm}$ is a  field of   Krein positive operators which is ``almost parallel'' in the sense that it satisfies:
%
%\be
%\Gamma_{v'}=H_{\cal P}(c_v^\ell(\Gamma_v))
%\ee
%
%for any path of length $\ell$ from $v$ to $v'$.
\end{theorem}

Clearly, the existence of a future-directed timelike and parallel vector field is equivalent to the existence of a timelike vector which is fixed by all the element of the holonomy group.

\begin{figure}[hbtp]
\begin{center}
 $$\xymatrix{ &3\ar[ddl]_c^{\rm Id}&    \cr
  & &\cr
1\ar[rr]^a_{\rm boost}  & & 2\ar[uul]_b^{\rm Id}}$$
\caption[toto]{Let us write the elements of ${\cal A}^0$ in the form $a=\pmatrix{a(1)&0\cr 0&a(2)}\oplus\pmatrix{a(2)&0\cr 0&a(3)}\oplus\pmatrix{a(3)&0\cr 0&a(1)}$. In view of theorem \ref{th6} a positive time-orientation form $\beta$ is completely determined by the Krein-positive operators $\Gamma_e^-$. Let us consider the simple case where $\Gamma_e^-=v_e$ is future-directed timelike vector. Then with the choices of discrete connection explained in the text we have  $a^{*_\beta}=\pmatrix{v_aa(1)^\times v_a^{-1}&0\cr 0&\Lambda(v_a)a(2)^\times\Lambda(v_a)^{-1}}\oplus\pmatrix{v_ba(2)^\times v_b^{-1}&0\cr 0&v_ba(3)^\times v_b^{-1}}\oplus\pmatrix{v_ca(3)^\times v_c^{-1}&0\cr 0&v_ca(1)^\times v_c^{-1}}$, which cannot be an element of ${\cal A}^0$ since we would then have  $v_a=v_b=v_c=\Lambda(v_a)$. Note that if we take $\Lambda$ to be a spatial rotation instead of a boost, this will work if we choose $v_a$ in the time axis fixed by this rotation.}
\label{counterexample}
\end{center}
\end{figure}
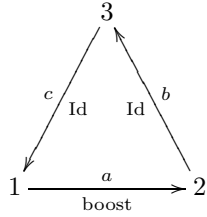

We can use the above theorem in particular to show the existence of non-reconstructible spectral spacetimes. The simplest we can imagine is defined over a triangle with edges $a,b,c$ and vertices $1,2,3$, and with $n=4$ (this example is detailed and illustrated in figure \ref{counterexample}). We take $g_1=g_2=g_3=\eta$, the Minkowski metric of signature $(1,3)$,   $S_1=S_2=S_3=S=\CC^4$, $h_b=h_c=Id_S$. To define $h_a$ we consider a boost $\Lambda$ and let $h_a$ be one of its two pre-images by the covering map $Spin^0(1,3)\rightarrow SO^+(1,3)$ (we identify $\CC l(\RR^{1,3})$ with $End(S)$ thanks to the representation $\rho$). Then  $\Lambda$  generates the holonomy group of $G$, and since it is a boost, it does not fix any timelike vector. 

% This example is illustrated in figure \ref{counterexample}.

% \newsavebox{\smlmat}% Box to store smallmatrix content
% \savebox{\smlmat}{$a=\pmatrix{a(1)&0\cr 0&a(2)}\oplus\pmatrix{a(2)&0\cr 0&a(3)}\oplus\pmatrix{a(3)&0\cr 0&a(1)}$}
% \newsavebox{\smlmat2}
% \savebox{\smlmat2}{$a^{*_\beta}=\pmatrix{v_aa(1)^\times v_a^{-1}&0\cr 0&\Lambda(v_a)a(2)^\times\Lambda(v_a)^{-1}}\oplus\pmatrix{v_ba(2)^\times v_b^{-1}&0\cr 0&v_ba(3)^\times v_b^{-1}}\oplus\pmatrix{v_ca(3)^\times v_c^{-1}&0\cr 0&v_ca(1)^\times v_c^{-1}}$}

% \begin{figure}[hbtp]
% \begin{center}
%  $$\xymatrix{ & 3 \cr
%  2 & \cr
%   & 1}$$
% \caption{Let us write the elements of ${\cal A}^0$ in the form ~\usebox{\smlmat}. In view of theorem \ref{theorem1} a positive time-orientation form $\beta$ is completely determined by the Krein-positive operators $\Gamma_e^-$. Let us consider the simple case where $\Gamma_e^-=v_e$ is future-directed timelike vector. Then with the choices of discrete connection explained in the text we have  ~\usebox{\smlmat2}, which cannot be an element of ${\cal A}^0$ since we would then have  $v_a=v_b=v_c=\Lambda(v_a)$. Note that if we take $\Lambda$ to be a spatial rotation instead of a boost, this will work if we choose $v_a$ in the time axis fixed by this rotation.}
% \label{counterexample}
% \end{center}
% \end{figure}

\subsubsection{Stable causality  when $n=4$}

Let us consider the issue of stable causality when $n=4$. We write $S_G^4$ for the split Dirac structure over a graph $G$, with the algebra ${\cal A}^0$ being the sections of the discrete bundle with fibre $\CC l(1,3)^0$. We assume here that $S_G^4$ satisfies the hypotheses of theorem \ref{th6}, and is thus a possibly non-reconstructible spectral spacetime.

Let us write $a\in{\cal A}^0$ in the form

\be
a(v)=f(v)1+\sum_{1\le i<j<4}g_{ij}(v)e_i(v)e_j(v)+h(v)\chi_v=f(v)+\alpha(v)+h(v)\chi_v\label{eq85}
\ee

where $f,g_{ij},h$ are scalar functions on $V$, $(e_i(v))_i$ is a pseudo-orthogonal basis in $X_v$, and the second equality defines $\alpha$. 

Let $\beta=[D,\pi(a)]$. We want $\beta^\times=\beta$ and $J\beta J^{-1}=-\beta$. Since

\bea
\beta^\times&=&-[D,\pi(\bar f)]-[D,\pi(\alpha)^\times]+[D,\pi(\bar h) \chi]\cr
J\beta J^{-1}&=&[D,\pi(\bar{f})]+[D,J\pi(\alpha)J^{-1}]-[D,\pi(\bar h)\chi]
\eea

we obtain that $[D,J\pi(\alpha)J^{-1}]=[D,\pi(\alpha)^\times]$. However, from (\ref{eq73}) we get

\bea
\pi(\alpha)^\times F(e,\pm)&=&h_e^\pm\alpha(e^\mp)^\times h_e^\mp F(e,\pm)\cr
&=&-h_e^\pm J_{e^\mp}\alpha(e^\mp)J_{e^\mp}^{-1}  h_e^\mp F(e,\pm),\mbox{ since }\alpha(e^\mp)^\times=-c(\alpha(e^\mp))\nonumber
\eea

and

\bea
J\pi(\alpha)J^{-1}F(e,\pm)&=&h_e^\pm J_{e^\mp}\alpha(e^\mp)J^{-1}F(e,\mp)\cr
&=&h_e^\pm J_{e^\mp}\alpha(e^\mp)J_{e^\mp}^{-1}h_e^\mp F(e,\pm)
\eea

hence $\pi(\alpha)^\times =-J\pi(\alpha)J^{-1}$. Thus $-[D,\pi(\alpha)^\times]=[D,J\pi(\alpha)J^{-1}]=[D,\pi(\alpha)^\times]=0$. Hence $[D,\pi(\alpha)]$ must vanish for $\beta$ to be Krein-self-adjoint and imaginary. We can then suppose without loss of generality that $\alpha=0$ when looking for a ``time function'' $a$.

% C'EST FAUX CAR \pi(a)+\pi(a)^\times n'est pas de la forme \pi(b) en général ! Thus we can suppose without loss of generality that $\pi(a)^\times=-\pi(a)$ and $J\pi(a)J^{-1}=-\pi(a)$ since the Krein-self-adjoint part as well as the $J$-real part of $\pi(a)$ commute  with $D$. Using (\ref{eq73}) we obtain

% \bea
% h_e^\pm a(e^\mp)^\times h_e^\mp&=&-a(e^\pm)\cr
% a(e^\mp)&=&-J_{e^\mp}h_e^\mp a(e^\pm)h_e^\pm J_{e^\mp}^{-1}
% \eea

Thanks to (\ref{eq61}) we see that the operators $\beta_e^\pm$ defined by (\ref{eq63}) for $\beta=[D,\pi(a)]$ are

\be 
\beta_e^\pm={1\over\delta_e}\big(\delta f_e\gamma_e^\pm   \pm \sigma h_e\chi_{e^\pm}\gamma_e^\pm \big)h_e^\pm:=\Gamma_e^\pm h_e^\pm
\ee

where 

\bea
\delta f_e&=&f(e^+)-f(e^-)\cr
\sigma h_e&=&h(e^+)+h(e^-)
\eea

These must be real functions in order for $\Gamma_e^\pm$ to be Krein-self-adjoint.

Let us write $\gamma_e^\pm$ in the form

\be
\gamma_e^\pm=v_e^\pm+\chi_{e^\pm}w_e^\pm
\ee
where $v_e^\pm,w_e^\pm\in X_{e^\pm}$ (real vectors). 

The operator $\Gamma_e^\pm$, which we ask to be Krein-positive, is thus

\be
\Gamma_e^\pm={1\over\delta_e}\left((\delta f_e v_e^\pm\pm\sigma h_e w_e^\pm)+\chi_{e^\pm}(\delta f_e w_e^\pm\pm \sigma h_e v_e^\pm)\right)
\ee

Thanks to lemma  \ref{lem8} in \cite{part1} we know that it is Krein-positive iff 

\bea
\delta f_e v_e^\pm\pm\sigma h_e w_e^\pm+(\delta f_e w_e^\pm\pm\sigma h_e v_e^\pm)&\in&C_{e^\pm}^+\cr
\delta f_e v_e^\pm\pm\sigma h_e w_e^\pm-(\delta f_e w_e^\pm\pm\sigma h_e v_e^\pm)&\in&C_{e^\pm}^+\label{eqpourrie}
\eea
where $C_{e^\pm}^+$ is  the forward light cone in $X_{e^+}^+$. These are four conditions, but thanks to (\ref{eq57}) we have

\bea
v_e^+&=&h_e^+v_e^-h_e^-\cr
w_e^+&=&-h_e^+w_e^-h_e^-\label{eq95}
\eea

hence (\ref{eqpourrie}) boils down to

% and (\ref{eq93}) translates as

% \bea
% (\delta f v_e^++\sigma h w_e^+)+\chi_{e^+}(\delta f w_e^++\sigma h v_e^+)&&\mbox{Krein }>0\cr
% (\delta f v_e^--\sigma h w_e^-)+\chi_{e^-}(\delta f w_e^--\sigma h v_e^-)&&\mbox{Krein }>0
% \eea

% Using (\ref{eq95}) we also obtain
% \bea
% (\delta f v_e^--\sigma h w_e^-)+\chi_{e^-}(-\delta f w_e^-+\sigma h v_e^-)&&\mbox{Krein }>0\cr
% (\delta f v_e^++\sigma h w_e^+)+\chi_{e^+}(-\delta f w_e^+-\sigma h v_e^+)&&\mbox{Krein }>0
% \eea

% Hence we obtain the two following conditions on $v_e^+,w_e^+$:
% \bea
% (\delta f v_e^++\sigma h w_e^+)+\chi_{e^+}(\delta f w_e^++\sigma h v_e^+)&&\mbox{Krein }>0\cr
% (\delta f v_e^++\sigma h w_e^+)-\chi_{e^+}(\delta f w_e^++\sigma h v_e^+)&&\mbox{Krein }>0
% \eea

% Once they are satisfied, the conditions on $v_e^-,w_e^-$ will automatically hold thanks to (\ref{eq95}). 

% Using the Dirac representation we easily see that both these conditions are equivalent to

\bea
(\delta f +\sigma h)(v_e^++w_e^+)&\in& C_{e^+}^+\cr
(\delta f -\sigma h)(v_e^+-w_e^+)&\in& C_{e^+}^+
\eea

Let us make the hypothesis that $v_e^+\pm w_e^+\in C_{e^+}^+$. Then the same will be true of $v_e^-\pm w_e^-$ by (\ref{eq95}). In this case we will say that the edge $e$ is timelike and future-directed. On such an edge we see that $f$ and $h$ must satisfy $\delta f>|\sigma h|$. We similarly define past-directed timelike edges. If $v_e^++w_e^+\in C_{e^+}^+$ and $v_e^+-w_e^+\in C_{e^+}^-$ we obtain that $\sigma h>|\delta f|$. We say that such an edge is of type $+\sigma$ . Finally if $v_e^++w_e^+\in C_{e^+}^-$ and $v_e^+-w_e^+\in C_{e^+}^+$ we have $\sigma h<-|\delta f|$ and we say that the edge is of type $-\sigma$.

 Here are some easy  conclusions we can draw from the above study. A \emph{timelike loop} is a loop of timelike edges which are all of the same future/past type.

\begin{propo}
\begin{itemize}

\item For $S_G^4$ to be stably causal it is necessary that $G$ does not contain any timelike loop.
 \item If all the edges of $G$ are timelike, it is also sufficient that $G$ does not contain any timelike loop. If all the edges of $G$ are of type $+\sigma$ (resp. $-\sigma$), then $S_G^4$ is always stably causal.
\item If $S_G^4$ is vectorial, it is stably causal iff it does not contain any timelike loop and all its edges are timelike.
\end{itemize}
\end{propo}
\begin{demo}
It is clear that $G$ cannot contain any timelike loop, since $f$ would have to satisfy $\delta f>|\sigma h|$ (if the loop consists of future-directed edges, $\delta f<-|\sigma h|$ if it consists of past-directed edges) on every edge of the loop, which is absurd 

The second assertion is immediate: it all the edges are timelike, we just give $G$ the orientation induced by the past/future type of its edges. Then we know that if it does not contain oriented loops for this orientation, we can find a strictly increasing function $f$ for the induced ordering, and we set $h$ to $0$. If all the edges are of $+\sigma$ (resp. $-\sigma$)-type we take $f=0$ and $h$ a positive (resp. negative) constant.

Finally if $S_G^4$ is vectorial we must have $(\delta f_e +\sigma h_e) v_e^+\in C_{e^+}^+$ and $(\delta f_e-\sigma h_e)v_e^+\in C_{e^+}^+$, hence $\delta f_e v_e^+\in C_{e^+}^+$. This proves that all edges are timelike and that $f$ must increase along future-directed timelike edges and decrease along past-directed ones. In other words $f$ must be a time function and this exists iff $G$ does not contain timelike loops. Conversely if there is no timelike loop and all the edges are timelike, we take $f$ to be a time function and set $h$ to $0$.
\end{demo}

We note that all the edges are timelike if 1) all the vectors $v_e$ are, and 2) the $w_e$ are small enough perturbations for $v_e\pm w_e$ to remain in the same half-cone as $v_e$. In such a case, the time function $a=f+h\chi$ also appears as a sum of a ``classical'' time function $f$ and a perturbation $h$ such that $\sigma h$ is small enough not to change the sign of $\delta f$. It remains to see if this nice characterization can be generalized to $n>4$.

On the contrary, the case  of graphs with edges of mixed types seems to be  much more intricate. The reader can check easily that all graphs with $3$ edges or less are stably causal, except the two timelike triangles. However the following graph of mixed type with $4$ edges and no timelike loop is not stably causal.  The graph is not drawn with the original orientation, which plays no role, but with the following convention:  a solid arrow which goes from bottom to top denotes a  future-directed timelike  edge  and  a  $+\sigma$ edge is dotted and goes from left to right.

 $$
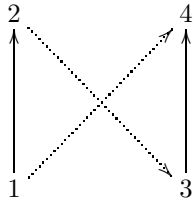
\xymatrix{
2\ar@{.>}[ddrr]& & 4\cr
&& \cr
1\ar[uu]\ar@{.>}[uurr]& &3\ar[uu]}$$

%$$\xymatrix{2\ar@{.>}[r] & 3\ar[d]  \cr
% 1\ar[u]\ar@{.>}[r]&4
%}$$

To see that this is not stably causal  simply note that the total variation $\Delta f$ of $f$ around the loop 1-2-3-4 must satisfy
$$0=\Delta f>|h(1)+h(2)|-(h(2)+h(3))+|h(3)+h(4)|-(h(1)+h(4))>0$$

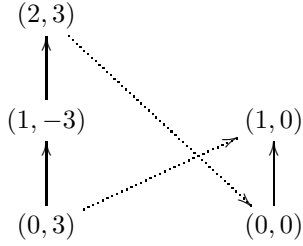
\begin{figure}[hbtp]
 $$\xymatrix{
{(2,3)}\ar@{.>}[ddrr]& & \cr
{(1,-3)}\ar[u]&&(1,0) \cr
{(0,3)}\ar[u]\ar@{.>}[urr]& &(0,0)\ar[u]}$$
\caption{In this graph we have indicated the value of a couple of functions $(f,h)$ at each vertex, which shows that $S_G^4$ is stably causal.}\label{figsc} 
\end{figure}

One could be under the impression that the problem comes from $+\sigma$ arrows joining vertices  which lie at ``different heights''. Indeed, if the graph contains only timelike or $+\sigma$ edges and no timelike loop, we can, neglecting the edges of type $+\sigma$, define the height $h(v)$ of a vertex $v$ to be the number of its predecessors. Then if the $+\sigma$ edges only join vertices which are at the same height, we can easily show that the couple $(f,h)$ with $f(v)=h(v)(h(v)+1)$ satisfies all the constraints. However the graph in figure \ref{figsc} shows the existence of other possibilities. 

We conclude here these preliminary explorations, leaving the question of stable causality of split Dirac structures essentially as an open problem.
%Since we are unable to pinpoint a  general principle, we will content ourselves by saying that vectorial graphs are ``more natural'' \ldots

\subsubsection{Relation with the discretization of the usual Dirac operator}\label{secdiscr}

We start by recalling the expression of the discretized Dirac operator introduced by Marcolli and Van Suijlekom in \cite{MVS}. Here it is supposed that $G$ is a graph embedded in a \emph{Riemannian} spin manifold $M$. For each vertex $v$, the space $X_v$ is $T_v M$, and $S_v$ is an irreducible representation space of $\CC l(T_v M,g_v)$. 

The discretized Dirac operator acts on the space $\Gamma (S)$ of sections  of the discrete spinor bundle $S=\bigsqcup_{v\in V}S_v$. The scalar  product of two elements of $\Gamma(S)$ is

\be
(\phi,\psi)=\sum_{v\in V}(\phi(v),\psi(v))\label{diagkrein}
\ee

and the discretized Dirac operator on $G$, which we write $\tilde D$, is

\be
(\tilde D\psi)(v)=i\sum_{e|t(e)=v}{1\over 2l_e}\gamma_e {\rm Hol}(e,\nabla^S)\psi(s(e))+i\sum_{e|s(e)=v}{1\over 2l_{\bar e}}\gamma_{\bar e}{\rm Hol}(\bar e,\nabla^S)\psi(t(e))\label{mvsdirac}
\ee

where  ${\rm Hol}(e,\nabla^S)$ is the parallel transport operator of the spin connection,  $l_e$ is the geodesic length of the embedded edge $e$, $\bar e$ is   $e$ with reversed orientation, and the gamma matrices are defined in the following way\footnote{In the original paper a factor of $i$ is included in the definition. We have left it outside to ease the comparison with our previous notations.}: at each vertex $v$ the tangent vectors $\dot e_i$ to the ingoing edges form a basis\footnote{In fact it can be an overcomplete system, but we discard this possibility here for simplicity's sake.} of $T_vM$ (which need not be orthonormal). We let $(\theta^i)$ be the corresponding dual basis of $T_vM^*$   and $\gamma_{e_i}=\rho_v((\theta^i)^\sharp)$ where $\rho_v : \CC l(T_v M,g_v)\rightarrow End(S_v)$ is the representation map.  Of course if $(\dot e_i)$ is an orthonormal basis, $\gamma_{e_i}$ is just $\rho_v(\dot e_i)$.

In the pseudo-Riemannian setting, we have to change the meaning of $l_e$ for the expression (\ref{mvsdirac}) to make sense. We  constrain the edge of the graph $G$ to be timelike or spacelike curves, and $l_e$ to be the proper time or proper distance along the edge, according to its nature. Let us observe once more that we cannot just take the discretized Dirac operator $\tilde D$ and build a spectral spacetime with it directly on $G$, using $\Gamma(S)$ with the diagonal Krein product (\ref{diagkrein}). Indeed, the 1-forms $\beta$ will have the same support as $ \tilde D$, in particular they have vanishing diagonal. Since $(.,.)$ is diagonal if one takes a spinor field  $\psi_v$ concentrated at $v$, $(\psi_v,\beta \psi_v)=0$, hence $\beta$ cannot turn $(.,.)$ into a scalar product.

Hence we must turn to  the split Dirac structure. The Dirac operator $D$ defined by (\ref{Dirac49}) acts on sections of the discrete spinor bundle $K$ over the split graph of $G$. To see the relation between $D$ and $\tilde D$ we introduce two natural maps between $\Gamma(S)$ and $K$.

First there is a natural  embedding of vector spaces $i : \Gamma(S)\rightarrow K$ defined by

$$i(\phi)(e,+)=\phi(t(e)),\qquad i(\phi(e,-))=\phi(s(e))$$

Note that this embedding does not respect the Krein product: if we write $\delta_v^{s_i}$ the element of $\Gamma(S)$ such that $\delta_v^{s_i}(v')=\delta_{v,v'}s_i$, and if we choose $(s_i)_{1\le i\le k}$ to be a pseudo-orthonormal basis of $S$, then the $\delta_v\otimes s_i$ will form a pseudo-orthonormal basis of $\Gamma(S)$. However an easy computation shows that 

$$(i(\delta_{v}^s),i(\delta_{v'}^{s'}))=\left\{\matrix{0,&\mbox{ if }v\mbox{ and }v'\mbox{ are not connected}\cr
(s,h_e^-s'),&\mbox{ if }(v,v')=e\cr
(s,h_e^+s'),&\mbox{ if }(v',v)=e}\right.
$$
 
We will call \emph{graph states} the elements of $i(\Gamma(S))$, since they are spinor fields on the split graph which can be defined directly on the original graph. We could hope that the Dirac operator $D$ defined by (\ref{Dirac49}) restrict to an operator on the graph states, but this is not the case: the image of a graph state by $D$ can fail to be a graph state. However we can also define the surjection $\Pi : K\rightarrow V$ which takes a function defined on $E\times \{+,-\}$ and averages over all the edges leading to, or coming from, a given vertex :

$$(\Pi F)(v)={1\over d_v}\left(\sum_{e|t(e)=v}F(e,+)+\sum_{e|s(e)=v}F(e,-)\right)$$

where $d_v$ is the degree of the vertex $v$. Since $\Pi\circ i=\id_{\Gamma(S)}$,  $i\circ \Pi$  is a projector with range  $i(\Gamma(S))$.  It turns out that we can adjust the scalars $\delta_e$ and the operators $\gamma_e^\pm$ and $h_e^\pm$ in (\ref{Dirac49}) so that  the following diagram commutes:

\be
\xymatrix{K\ar[r]^D&K\ar[d]^\Pi\cr
\Gamma(S)\ar[u]^i\ar[r]^{k\tilde D}&\Gamma(S)\label{eq101}
}
\ee

where $k$ is some constant. Indeed, we have:

\bea
(\Pi\circ D\circ i)(\psi)(v)&=&{1\over d_v}\big(\sum_{e|t(e)=v}D i(\psi)(e,+)+\sum_{e|s(e)=v}D i(\psi)(e,-)\big)\cr
&=&{1\over d_v}\big(\sum_{e|t(e)=v}{-1\over \delta_e}\gamma_e^+h_e^+ i(\psi)(e,-)+\sum_{e|s(e)=v}{1\over \delta_e}\gamma_e^-h_e^- i(\psi)(e,+)\big)\cr
&=&{1\over  d_v}\big(-\sum_{e|t(e)=v}{1\over \delta_e}\gamma_e^+h_e^+ \psi(s(e))+\sum_{e|s(e)=v}{1\over \delta_e}\gamma_e^-h_e^- \psi(t(e))\big)\nonumber
\eea

We see that   this expression is essentially the same as (\ref{mvsdirac}). Indeed,    in our definition of $D$ we can choose $\gamma_e^+$ and $\gamma_e^-$ to be $\gamma_e$ and $-\gamma_{\bar e}$ as defined by Marcolli and van Suijlekom,  $\delta_e$ to be the proper time/distance along $e$, and $h_e^\pm$ to be the parallel tranport operators of the spin connection along $e/\bar e$, which are inverse of each other and Krein unitary. Then the two expressions will coincide apart from the numerical factor $-i{ d_v\over 2}$. This factor will not depend on $v$ if we ask the gamma matrices to form a basis of spacetime at each vertex, as we have done.
 
Hence diagram (\ref{eq101}) commutes up to a constant. But there is more: since we have chosen $h_e$ to be ${\rm Hol}(e,\nabla^S)$, the discrete connection will automatically be metric, spin,  and orientation preserving. Moreover, the choice of the $\gamma_e^\pm$ puts us in the complete vectorial case. Hence the discrete connection is spin. Condition (\ref{eq83}) is the only one remaining to ensure that the hypothesis of theorem \ref{th6} are satisfied, and here it takes the form 

\be
\gamma_{e}=-{\rm Hol}(e,\nabla^S)(\gamma_{\bar e})\label{eq103}
\ee

If we ask the tangent vectors $\dot e_i$ to form a pseudo-orthonormal basis, we see, taking into account the fact that $e$ and $ \bar e$ have opposite orientation, that (\ref{eq103})  is equivalent to the tangent vector to $e$ at $t(e)$ being the parallel tranport of the tangent vector to $e$ at $s(e)$. This will be automatically fulfilled if we furthermore ask $e$ to be a geodesic segment. We then arrive at the following conclusion: \emph{when the graph $G$ is made of  ``pseudo-orthonormal'' geodesic segments, the choices we need to make (\ref{eq101}) commute (up to the factor $-in/2$) automatically turn the split Dirac structure into an antilorentzian spectral spacetime.} This will be reconstructible or not according to the properties of the curvature of the spacetime being approximated.

\section{Conclusion and outlook}\label{conclusion}

In this paper we have given a new definition extending spectral triples to the Lorentzian/antilorentzian cases. Actually, the novelty lies only in these few points: we do not consider a Hilbert space to be a basic object in the structure, or, equivalently, we do not fix a particular time-orientation form, and   we do not impose the commutation of this time orientation form with the algebra. Moreover we do not suppose that the time-orienation forms need to be normalized, hence becoming fundamental symmetries, and we ask that they anticommute with the charge conjugation. These two last points are of   lesser importance, but the the first two turn out to have a dramatic consequence from a mathematical point of view: the algebra under consideration is generally neither a $C^*$ nor a Krein $C^*$-algebra. On the other hand, we think that the need for a  change of focus from Hilbert spaces to Krein spaces is better appreciated from the point of view of physics. If we draw an analogy with special relativity, the Krein space plays the role of Minkowski space, and the time-orientation forms stand for the inertial reference frames\footnote{Of course this is just an analogy: in our case the time-orientation forms are associated to congruence of timelike curves, i.e. partitions of spacetime into observers worldline, and these observers need not be locally inertial.}. One can analyze the relativistic phenomena without making use of reference frames at all, but if we want to use them we will need at least two, so it would seem to be a very bad idea to fix a particular one in the background. 

In fact we find it remarkable that, even though our argumentation has been essentially mathematical,   the end product is so consonant with physics. Indeed, already in Euclidean noncommutative geometry the  kind of manifolds which are generalized to the noncommutative settings are spin manifolds, i.e. those on which matter fields can exist. This   may be an underestimated sign that noncommutative geometry is on the right track. In the Lorentzian case we have seen that the class of manifolds is further restricted to time and space orientable ones, i.e. manifolds on which chirality and the arrow of time can be  defined. We think this is also a good sign, especially because the time-orientations are associated with the Hilbert structure and most interpretations of quantum mechanics depend on a time direction. %[CUT : this is a bit cryptic... rephrase or cut]

However, we cannot content ourselves with omens: the only way we can be sure the approach is physically correct is by using to build a consistent physical model. We will certainly need to extend the work presented here to other signatures, since there are indications that the total signature one needs for the spectral standard model is neither Euclidean nor Lorentzian \cite{vdd}, \cite{nadir}. The calculations of the spectral action and its variation also need to be done. On this last point we have already performed such a calculation on a vectorial split Dirac structure, and the results are encouraging, but they need to be generalized.

On a more mathematical side, a lot remains to be done also. Let us list a few problems that we have left open, in no particular order:

\begin{enumerate}
\item What are the correct analytical requirements to be put on the spectral spacetimes ? 
\item How can we extend the approach to cover the odd dimensional case ?
\item Are all time-orientation forms normalizable ?
\item Is it possible to characterize stably causal split Dirac structure   ?
\item The finite-dimensional  spectral triples can be classified (\cite{kraj}, \cite{PS98}, \cite{cac}). How does this classification generalize to finite-dimensional spectral spacetimes ?
\item Is it possible to extend the definition of spectral spacetime to other signatures, taking care of junk forms ? 
\item Does there exists a structure satisfying all the properties of a spectral spacetime but with ``wrong KO-signs'' ?
\end{enumerate}

To these questions we might add a particularly important one, in our opinion: can one use the split Dirac structure to build a model of quantum spacetime ? If indeed spacetime is discrete at a quantum level, and can be analyzed with the tools of noncommutative geometry, then we must take into account the fact that a discrete approximation of a spacetime in noncommutative geometry is not given by a finite graph embedded in the manifold, but uses a split graph instead. %This seems to us an inescapable conclusion that split graphs must be used in approaches like... CITER VS + M

Finally, it has to be said that the former approach of noncommutative causality initiated by the author (see \cite{bes0}) using isocones has to be updated to take into account the shift of emphasis from Hilbert to Krein space that we have advocated here.

% Our argumentation has been essentially mathematical, but we find it remarkable that the end product fits so nicely with physics: on the one hand the kind of manifolds which are generalized to the noncommutative settings are those on which matter fields can exist and on which chirality as well as the distinction of the future and past have a meaning. On the other hand the Hilbert structure is seen to be relative to classes of observers, expressed by the time-orientation 1-forms. This is a lesson which can also be drawn from relativistic quantum mechanics. 

% Important : Hilbert space $\rightarrow $ positivity $\rightarrow $ probability.

% Hilbert space = Krein space + Krein positive operator. Krein positive operator = observers. Without observers, no probability ? No quantum mechanics ?

\end{document}